\documentclass{gtart_h}  

\def\ifplaintex{\expandafter\ifx\csname documentclass\endcsname\relax}

\def\ifplaintex{\expandafter\ifx\csname documentclass\endcsname\relax}

\ifplaintex 
\else
\fi

\expandafter\ifx\csname epsfbox\endcsname\relax\input epsf\fi

\def\gt{{\mathsurround=0pt\it $\cal G\mskip-2mu$eometry \&\ 
$\cal T\!\!$opology}}        

\def\gtp{{\mathsurround=0pt\it $\cal G\mskip-2mu$eometry \&\ 
$\cal T\!\!$opology $\cal P\!$ublications}}  

\def\lognumber#1{\def\thelognumber{#1}}
\def\volumenumber#1{\def\thevolumenumber{#1}}
\def\papernumber#1{\def\thepapernumber{#1}}
\def\volumeyear#1{\def\thevolumeyear{#1}}

\def\pagenumbers#1#2{\def\startpage{#1}\def\finishpage{#2}}
\def\published#1{\def\publishdate{#1}}
\def\proposed#1{\def\theproposer{#1}}
\def\seconded#1{\def\theseconders{#1}}
\def\received#1{\def\receiveddate{#1}}
\def\revised#1{\def\reviseddate{#1}}
\def\accepted#1{\def\accepteddate{#1}}
\def\asciititle#1{\def\theasciititle{#1}}

\def\asciiaddress#1{\def\theasciiaddress{#1}}
\def\asciiemail#1{\def\theasciiemail{#1}}

\long\def\asciiabstract#1{\long\def\theasciiabstract{#1}}
\def\asciikeywords#1{\def\theasciikeywords{#1}}
\def\shortauthors#1{\def\theshortauthors{#1}}
\def\shorttitle#1{\def\theshorttitle{#1}}

\let\\\par\let\thelognumber\relax
\let\thevolumenumber\relax\let\thepapernumber\relax
\let\thevolumeyear\relax\let\thesamplenumber\relax\let\startpage\relax
\let\finishpage\relax\let\publishdate\relax\let\receiveddate\relax
\let\reviseddate\relax\let\accepteddate\relax\let\theasciititle\relax
\let\theasciiauthors\relax\let\theasciiaddress\relax
\let\theasciiabstract\relax\let\theasciikeywords\relax
\let\theasciiemail\relax\let\theshortauthors\relax\let\theshorttitle\relax

\long\def\maketitlep{   

\count0=\startpage

\gt\hfill      
\hbox to 77pt{\vbox to 0pt{\vglue -15pt\epsfbox{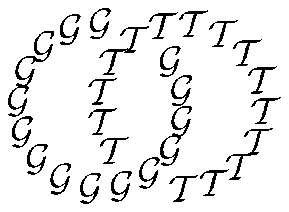}\vss}\hss}
\break
{\small\ifx\thesamplenumber\relax 
Volume \else Sample
\fi\thevolumenumber\ (\thevolumeyear)
\startpage--\finishpage\nl
Published: \publishdate}
\vglue 0.5truein plus 0.4fil minus 0.1truein

{\parskip=0pt\leftskip 0pt plus 1fil\def\\{\par\smallskip}{\ifplaintex\large
\else\Large\fi\bf\thetitle}\par\medskip}   

\vglue 0pt plus 0.1fil 

{\parskip=0pt\leftskip 0pt plus 1fil\def\\{\par}{\sc\theauthors}
\par\medskip}

\vglue 0pt plus 0.1fil 

{\small\parskip=0pt\let\newline\\
{\leftskip 0pt plus 1fil\def\\{\par}{\sl\theaddress}\par}
\expandafter\ifx\theemail\relax    
\relax\else\vglue 5pt plus 0.02fil minus 2pt\def\\{\stdspace{\rm 
and}\stdspace} 
\cl{Email:\stdspace\tt\theemail}\fi
\ifx\theurl\relax                  
\relax\else\vglue 5pt plus 0.02fil minus 2pt\def\\{\stdspace{\rm 
and}\stdspace}
\cl{URL:\stdspace\tt\theurl}\fi\par}

\vglue 7pt plus 0.3fil minus 3pt

{\bf Abstract}
\vglue 5pt plus 0.1fil minus 2pt

\theabstract

\vglue 7pt plus 0.3fil minus 3pt

{\bf AMS Classification numbers}\quad Primary:\quad \theprimaryclass

Secondary:\quad \thesecondaryclass

\vglue 5pt plus 0.3fil minus 2pt

{\bf Keywords:}\quad \thekeywords

\vglue 10pt plus 0.5fil minus 5pt

{\small  Proposed: \theproposer\hfill Received: \receiveddate\nl
Seconded: \theseconders\hfill 
\ifx\reviseddate\relax                         
Accepted: \accepteddate                        
\else
Revised: \reviseddate                          
\fi}
\eject
}       

\font\phead=cmsl9 scaled 950
\font=cmsl9 scaled 1050
\font\pnum=cmbx10 scaled 913
\font\lnum=cmbx10 
\font\pfoot=cmsl9 scaled 950
\font=cmsl9 scaled 1050
\ifplaintex
\headline{\vbox to 0pt{\vskip -4.5mm\line{\small\phead\ifnum
\count0=\startpage ISSN 1364-0380 (on line)
1465-3060 (printed) \hfill {\pnum\folio}\else\ifodd\count0\def\\{ }%
\ifx\theshorttitle\relax\thetitle\else\theshorttitle\fi\hfill{\pnum\folio}
\else\def\\{ and }{\pnum\folio}\hfill\ifx\theshortauthors\relax\theauthors
\else\theshortauthors\fi\fi\fi}\vss}}
\footline{\vbox to 0pt{\vglue 0mm\line{\small\pfoot\ifnum\count0=\startpage
\copyright\ \gtp\hfill\else
\gt, Volume \thevolumenumber\ (\thevolumeyear)\hfill\fi}\vss
}}
\else
\makeatletter
\def\@oddhead{{\small\lhead\ifnum\count0=\startpage ISSN 1364-0380 (on line)
1465-3060 (printed) \hfill {\lnum\number\count0}\else\ifodd\count0
\def\\{ }\ifx\theshorttitle\relax \thetitle \else\theshorttitle\fi\hfill
{\lnum\number\count0}\else\def\\{ and }{\lnum\number\count0}
\hfill\ifx\theshortauthors\relax 
\theauthors\else\theshortauthors\fi\fi\fi}}\def\@evenhead{\@oddhead}
\def\@oddfoot{\small\lfoot\ifnum\count0=\startpage\copyright\ \gtp\hfill\else
\gt, Volume \thevolumenumber\ (\thevolumeyear)\hfill\fi}
\def\@evenfoot{\@oddfoot}
\makeatother
\fi

\newwrite\gtoutfile
\long\gdef\makeheadfile{  
{\def\\{, }\def\s{ }
\immediate\openout\gtoutfile head.xxx
\immediate\write\gtoutfile{Proxy-for: \ifx\theasciiauthors\relax
\theauthors\else\theasciiauthors\fi\s<\ifx\theasciiemail\relax\theemail\else\theasciiemail\fi>}
\immediate\write\gtoutfile{\noexpand\\}
\immediate\write\gtoutfile{Authors: \ifx\theasciiauthors\relax
\theauthors\else\theasciiauthors\fi}
{\def\\{ }\immediate\write\gtoutfile{Title: \ifx\theasciititle\relax
\thetitle\else\theasciititle\fi}}
\immediate\write\gtoutfile{Subj-class: GT or SG or MG etc}
\immediate\write\gtoutfile{MSC-class: \theprimaryclass\ifx\thesecondaryclass\relax\else, \thesecondaryclass\fi}
\immediate\write\gtoutfile{Journal-ref: Geom. Topol. \thevolumenumber
(\thevolumeyear) \startpage-\finishpage}
\immediate\write\gtoutfile{Comments: Published by Geometry and Topology at}
\immediate\write\gtoutfile{\s\s http://www.maths.warwick.ac.uk/gt/GTVol\thevolumenumber/paper\thepapernumber.abs.html}
\immediate\write\gtoutfile{\noexpand\\}
\immediate\write\gtoutfile{}
\ifx\theasciiabstract\relax
\immediate\write\gtoutfile{\theabstract}\else
\immediate\write\gtoutfile{\theasciiabstract}\fi
\immediate\write\gtoutfile{}
\immediate\write\gtoutfile{\noexpand\\}
\immediate\write\gtoutfile{}
\immediate\closeout\gtoutfile}}  

\def\maketitlepage{\maketitlep\makeheadfile}
\let\maketitle\maketitlepage

\lognumber{249}
\received{24 April 2002}
\volumenumber{9}\papernumber{27}\volumeyear{2005}
\pagenumbers{1187}{1220}   
\revised{28 April 2005}
\published{1 July 2005}
\accepted{6 June 2005}
\proposed{Ralph Cohen}
\seconded{Gunnar Carlsson, Haynes Miller}

\usepackage{amssymb,amsmath}
\usepackage{enumerate}
\usepackage[all]{xy}
\input diagxy
\newtheorem{theorem}{Theorem}[section]
\newtheorem{lemma}[theorem]{Lemma}
\newtheorem{proposition}[theorem]{Proposition}
\newtheorem{corollary}[theorem]{Corollary}

\theoremstyle{definition}

\newtheorem{definition}[theorem]{Definition}
\newtheorem{example}[theorem]{Example}

\theoremstyle{remark}

\newtheorem*{example*}{Example}
\newtheorem{note}[theorem]{Note}
\newtheorem{inote}[theorem]{Interesting note}

\newcommand{\calc}{\mathcal C}
\newcommand{\calz}{\mathcal Z}
\newcommand{\cyc}[2]{{\calz}^{#1}(#2)}
\newcommand{\cqt}[2]{{\calz}_{\bbh}^{#1}{(#2)}}
\newcommand{\cqtav}[2]{{\calz}^{#1}{(#2)}^\text{av}}
\newcommand{\cqtrd}[2]{\widetilde{\calz}_{\bbh}^{#1}{(#2)}}
\newcommand{\csus}{\Sigma \kern-.6em\raise1pt\hbox{/}}
\newcommand{\chv}[3]{{\calc}^{#1}_{#2}(#3)}
\newcommand{\jj}{\mathbf j}
\newcommand{\quat}{\gamma}
\newcommand{\quatt}{\Gamma}
\newcommand{\real}{\phi}
\newcommand{\reall}{\Phi}
\newcommand{\qa}{quaternionic\ }
\newcommand{\Qa}{Quaternionic\ }

\newcommand{\gB}{\mathfrak B}
\newcommand{\gR}{\mathfrak R}

\newcommand{\bbz}{\mathbb Z}
\newcommand{\bbr}{\mathbb R}
\newcommand{\bbc}{\mathbb C}
\newcommand{\bbh}{\mathbb H}
\newcommand{\bbp}{\mathbb P}
\newcommand{\bbd}{\mathbb D}
\newcommand{\cu}{\mathcal U}
\newcommand{\cc}{\mathcal C}
\newcommand{\call}{\mathcal L}
\newcommand{\cf}{\mathcal F}
\newcommand{\cz}{\mathcal Z}

\newcommand{\cco}[1]{{\mathcal O}_{#1}}

\newcommand{\ct}{\mathcal T}
\newcommand{\ci}{\mathcal I}

\renewcommand{\b}{\beta}
\newcommand{\g}{\gamma}
\newcommand{\la}{\lambda}
\newcommand{\G}{\Gamma}
\newcommand{\naive}{na{\"\i}ve }

\newcommand{\equdef}{\overset{\text{def}}{=}}

\newcommand{\codim}{\operatorname{codim}}

\newcommand{\arr}{\longrightarrow}
\newcommand{\Hom}{\operatorname{Hom}}
\newcommand{\Mor}{\operatorname{Mor}}

\newcommand{\Div}{\operatorname{Div}}

\newcommand{\sm}{\operatorname{Sym}^2}
\newcommand{\pch}[1]{\bbp_{\bbc}(\bbh^{#1})}
\newcommand{\zt}{$\bbz_2$--}
\newcommand{\pt}{\widetilde{\pi}}

\newcommand{\lo}{\ell_0}
\newcommand{\BP}{{\mathbb P}}
\newcommand{\BD}{{\mathbb D}}
\newcommand{\BC}{{\mathbb C}}
\newcommand{\BH}{{\mathbb H}}
\newcommand{\BZ}{{\mathbb Z}}
\newcommand{\tbbd}{\widetilde{\bbd}}

\newcommand{\plo}{p^*\lo}

\newcommand{\U}{{\mathcal Q}}

\newcommand{\ob}{\overline{\U}}

\newcommand{\aand}{\qquad\text{and}\qquad}
\newcommand{\ev}{\operatorname{ev}}
\newcommand{\odd}{\operatorname{odd}}
\newcommand{\grpi}{\gR}
\newcommand{\tF}{\widetilde{F}}
\newcommand{\hF}{\widehat F}

\newcommand{\Id}{\operatorname{Id}}
\newcommand{\Gal}{\operatorname{Gal}}
\renewcommand{\to}{\rightarrow}

\begin{document}

\title{Algebraic cycles and the classical groups\\Part II: Quaternionic cycles}
\asciititle{Algebraic cycles and the classical groups II: Quaternionic cycles}
\shorttitle{Algebraic cycles and the classical groups II: Quaternionic cycles}

\authors{H Blaine Lawson Jr\\Paulo Lima-Filho\\Marie-Louise Michelsohn}
\shortauthors{Lawson, Lima-Filho and Michelsohn}

\address{{\rm BL, MM:} Department of Mathematics, Stony Brook University\\Stony
Brook, NY 11794, USA\\\smallskip
\\{\rm PL:} Department of Mathematics, Texas A{\&}M
University\\College Station, TX 77843, 
USA\\\smallskip\\{\tt\mailto{blaine@math.sunysb.edu}, 
\mailto{plfilho@math.tamu.edu}, \mailto{mlm@math.sunysb.edu}}}

\asciiaddress{BL, MM: Department of Mathematics, Stony Brook University\\Stony
Brook, NY 11794, USA\\
\\PL: Department of Mathematics, Texas A\&M
University\\College Station, TX 77843, 
USA}

\asciiemail{blaine@math.sunysb.edu, plfilho@math.tamu.edu,
  mlm@math.sunysb.edu}

\begin{abstract} 
In part I of this work we studied the spaces of real algebraic 
cycles on a complex projective space  $\mathbb{P}(V)$, 
where $V$ carries a real structure, and completely determined their homotopy type.
We also extended some functors in $K$--theory to algebraic
cycles, establishing a direct relationship to characteristic
classes for the classical groups, specially Stiefel--Whitney
classes. In this sequel, we establish
corresponding results in the case where $V$ has a quaternionic
structure. The determination of the homotopy type of quaternionic
algebraic cycles is more involved than in the real case, but has a
similarly simple description. The stabilized space of quaternionic algebraic
cycles admits a nontrivial infinite loop space structure
yielding, in particular, a delooping of the total Pontrjagin class
map. This stabilized space is directly related to an extended notion of
quaternionic spaces and bundles ($KH$--theory), in analogy with
Atiyah's real spaces and $KR$--theory, and the characteristic classes
that we introduce for these objects are nontrivial. The paper ends
with various examples and applications.
 \end{abstract}

\asciiabstract{%
In part I of this work we studied the spaces of real algebraic cycles
on a complex projective space P(V), where V carries a real structure,
and completely determined their homotopy type.  We also extended some
functors in K-theory to algebraic cycles, establishing a direct
relationship to characteristic classes for the classical groups,
specially Stiefel-Whitney classes. In this sequel, we establish
corresponding results in the case where V has a quaternionic
structure. The determination of the homotopy type of quaternionic
algebraic cycles is more involved than in the real case, but has a
similarly simple description. The stabilized space of quaternionic
algebraic cycles admits a nontrivial infinite loop space structure
yielding, in particular, a delooping of the total Pontrjagin class
map. This stabilized space is directly related to an extended notion
of quaternionic spaces and bundles (KH-theory), in analogy with
Atiyah's real spaces and KR-theory, and the characteristic classes
that we introduce for these objects are nontrivial. The paper ends
with various examples and applications.}


\asciikeywords{Quaternionic algebraic cycles, characteristic classes,
 equivariant infinite loop spaces, quaternionic K-theory}

{\small\maketitlepage}

\section{Introduction}

In Part 1 \cite{QQLLM} of this paper  we studied the spaces of
real algebraic cycles in $\bbp(\bbc^n)$ and found that their
homotopy structure was particularly simple and surprisingly
related to the Stiefel--Whitney and Pontrjagin classes.  We saw
that the stabilized space $\cz^{\infty}_{\bbr}$ of all such cycles
is an $E_{\infty}$--ring space whose homotopy groups
$\pi_*\cz^{\infty}_{\bbr}$ form a graded ring isomorphic to
$\bbz[x,y]/(2y)$.  Furthermore, the standard complexification and
forgetful functors in $K$--theory were shown to extend over the
characteristic homomorphisms to infinite loop maps of cycle
spaces. Here in Part 2 we shall establish analogous results for
spaces of quaternionic cycles.

Recall that a real vector space is a pair $(V,\rho)$ where $V$ is a
complex vector space  and  $\rho\co V\to V$ is  an anti-linear map
with $\rho^2 = \Id$.  A real algebraic subvariety (or more
generally a real algebraic cycle) in $\bbp(V)$ is one which is
fixed by the induced involution   $\rho \co \bbp(V)\to \bbp(V)$ (the
action of $\Gal(\bbc/\bbr)$). This condition is equivalent to
assuming that the subvariety is defined by real algebraic
equations.

Analogously a quaternionic vector space is a pair $(V,\jj)$ where
$\jj\co V\to V$ is an anti-linear map with $\jj^2 =  -\Id$.  A
quaternionic algebraic variety (or cycle) in $\bbp(V)$ is one
which is fixed by the induced involution.  Such subvarieties are
distinctly different from real ones.  The map $\jj$ is fixed-point
free on $\bbp(V)$ and therefore on every subvariety.  The induced
antilinear bundle map on ${\mathcal O}(1)$ has square $-\Id$.

In \cite{QLLM} we prove that the quaternionic suspension of cycles
gives a $\bbz_2$--homotopy equivalence $ \csus_{\bbh}\co  
\cz^q(\bbp(V)) \to\cz^q(\bbp(V\oplus \bbh)) $
 of cycle spaces and in particular a homotopy equivalence
$$
\csus_{\bbh} \co  \cz^q_{\bbh}(\bbp(V))\arr
\cz^q_{\bbh}(\bbp(V\oplus \bbh))
$$
of the subgroups of codimension--$q$ quaternionic cycles. For $q$
odd one can thereby reduce to 0--cycles and apply Dold--Thom
\cite{QDT,QQQLi} to determine the  structure of these
spaces. However, for $q$ even one can only reduce to $1$--cycles by
this method. The determination of the homotopy type of
$\cz^q_{\bbh}$ when $q$ is even, is one of the main results of
this paper (Theorem \ref{thm:2.3}). Its proof, which involves new
constructions and techniques, is given in Section~\ref{sec:6}.

Under stabilization there are two   limits, $\cz_{\bbh}^{2\infty}$
and $\cz_{\bbh}^{2\infty+1}$, whose homotopy groups together form a
$\bbz_2{\times}\bbz$--graded ring under the algebraic join pairing.
Our second   result is the determination of this ring which turns
out to be quite simple. (See Theorem \ref{thm:3.4}.)

In analogy with the real and complex cases, we then show that the
inclusion of linear quaternionic cycles (the quaternionic
Grassmannian) into the space of all cycles yields the
characteristic map $\mathbf{BSp} \arr \prod_i K(\bbz,4i)$ classifying
the total Pontrjagin class. Furthermore, the forgetful functor and
the quaternionification functors from $K$--theory are shown to extend
over the characteristic maps to infinite loop maps of cycles
spaces.

It turns out that the spaces $\cz_{\bbh}^{2\infty}$ and
$\cz_{\bbh}^{2\infty+1}$ have a second, more mysterious
$\bbz_2$--component which is not seen by the characteristic map
from $\mathbf{BSp}$.  However, there is an extended notion of
quaternionic spaces and  bundles, in analogy with Atiyah's notion
of real spaces and real bundles \cite{QA}; and for such creatures
our new $\bbz_2$--characteristic classes are nontrivial. This is
discussed at the end of the paper where   examples and
applications are given.

\section{Spaces of quaternionic cycles}
\label{sec:2}

A \emph{quaternionic structure} on a complex vector space $V$ is a
$\bbc$--antilinear map $\jj \co  V \rightarrow V$ such that $\jj^2 =
-1.$   A \emph{quaternionic vector space} is a pair $(V,\jj)$
consisting of a complex   vector space $V$ and a quaternionic
structure $\jj$. Any quaternionic vector space is equivalent to
$(\bbh^n, \jj_0)$ where $\jj_0$ is left scalar multiplication by
the quaternion $j$.

A  quaternionic structure $\jj$  on $V$ induces a free
anti-holomorphic involution  $\jj\co \bbp(V) \to \bbp(V)$ which can
be viewed as follows.  Let $\pi\co \bbp(V) \to \bbp_{\bbh}(V)$ be the
projection from the complex to the quaternionic projective space
of $V$ whose fibres are projective lines. Then $\jj$ preserves the
fibres of $\pi$ and acts on them as the antipodal map on $S^2$.
This map $\jj$ induces an anti-holomorphic involution on the Chow
varieties $\chv{q}{d}{\bbp(V)}$, which in turn induces an
automorphism
\begin{equation}
\label{eq:2.1}
\jj \co  \cyc{q}{\bbp(V)} \rightarrow \cyc{q}{\bbp(V)}.
\end{equation}
of the topological group of all codimension--$q$ cycles on
$\bbp(V)$.

\begin{definition}
\label{def:2.1}
Let $(V,\jj)$ be a quaternionic vector space. Then the group
$\cqt{q}{\bbp(V)}$ of \emph{quaternionic algebraic cycles} of
codimension $q$ on $\bbp(V)$, is the fixed point set  of the
involution  \eqref{eq:2.1}. It contains the   closed subgroup of
\emph{averaged cycles} $\cqtav{q}{\bbp(V)}= \{ c + \jj c \; | \;
c\in \cyc{q}{\bbp(V)}\}$. We define the group of \emph{reduced
quaternionic algebraic cycles}
 to be the quotient
$$
\cqtrd{q}{\bbp(V)} = \cqt{q}{\bbp(V)} /    \cqtav{q}{\bbp(V)}.
$$
Note that $\cqtrd{q}{\bbp(V)}$ is the topological $\bbz_2$--vector
space freely generated by the $\jj$--invariant irreducible
subvarieties of codimension--$q$ in $\bbp(V)$.
\end{definition}

\begin{example*}
When $q = \dim(V)-2$, so that the cycle dimension is 1, then the
basis elements of $\cqtrd{q}{\bbp(V)}$ are irreducible,
$\jj$--invariant algebraic curves.  Note that if $C\subset \bbp(V)$
is such a curve, then the quotient $C_0 =  C/\bbz_2$ of $C$ by
$\jj$ is a ``nonorientable'' algebraic curve carrying a
``nonorientable'' conformal structure.  Conversely, any such
creature $C_0$ has a \zt covering $C\to C_0$ by a complex analytic
curve where the covering involution $\jj\co  C\to C$ is
anti-holomorphic and free. The natural embedding spaces for such
objects are $(\bbp(V), \jj)$ where $V$ is quaternionic. Thus one
could think of $\cqtrd{q}{\bbp(\bbh^2)}$ as the topological \zt
vector space generated by the irreducible ``nonorientable
algebraic curves'' in  $\bbp^3$.
\end{example*}

Given a quaternionic vector space $(V, \jj)$ and an algebraic
subvariety $Z\subset \bbp(V)$, we define the \emph{quaternionic
algebraic suspension} $\csus_{\bbh}(Z) \subset \bbp(V\oplus \bbh)$
to be the union of all complex projective lines joining $Z$ to
$\bbp(\{0\}\oplus \bbh)$.  This determines a continuous \zt
equivariant homomorphism (cf \cite[Section~6]{QLLM})
\begin{equation}
\label{eq:2.2}
\csus_{\bbh} \co   \cyc{q}{\bbp(V)}  \arr  \cyc{q}{\bbp(V\oplus
\bbh)}.
\end{equation}

\begin{theorem}{\rm \cite[Theorem~6.1]{QLLM}}\qua
\label{thm:2.2}
The quaternionic algebraic suspension homomorphism \eqref{eq:2.2}
is a \zt homotopy equivalence.  It induces homotopy equivalences
$$\csus_{\bbh} \co   \cqt{q}{\bbp(V)} \xrightarrow{\ \cong\ }
\cqt{q}{\bbp(V\oplus \bbh)}$$
$$\csus_{\bbh} \co  \cqtrd{q}{\bbp(V)} \xrightarrow{\ \cong\ }
\cqtrd{q}{\bbp(V\oplus \bbh)}\leqno{\hbox{and}}
$$
for all $q < \dim_{\bbc}(V)$.
\end{theorem}

Thus the homotopy types of the spaces $\cqt{q}{\bbp(V)}$ and
$\cqtrd{q}{\bbp(V)}$  depend only on $q$.  One of our main results
is the following computation of these homotopy types. The proof is
given in Section \ref{sec:6}.

\begin{theorem}
\label{thm:2.3}
For any quaternionic vector space $(V,\rho)$ there are canonical
homotopy equivalences:
\begin{enumerate}
\item[{\rm(i)}] $ \cz^{2q}_{\bbh}(\bbp(V)) \cong \prod_{j=0}^{q} K(\bbz,
4j) \times \prod_{j=1}^{q} K(\bbz_2, 4j{-}2)$
and $\cqtrd{2q}{\bbp(V)} \cong \bbz_2$,
\item[{\rm(ii)}] $ \cz^{2q+1}_{\bbh}(\bbp(V)) \cong \prod_{j=0}^{q} K(\bbz,
4j) \times \prod_{j=0}^{q} K(\bbz_2, 4j{+}1)$
and $\cqtrd{2q+1}{\bbp(V)} \cong \{\text{point}\} $ for all $q$.
\end{enumerate}
\end{theorem}

Note that $\csus_{\bbh}$ changes cycle dimension by 2.  Therefore,
for cycles of odd codimension we can de-suspend to the case of
0--cycles and apply the Dold--Thom Theorem.  For cycles of even
codimension one must take a complex suspension by ${\mathcal
O}(2)$ and then de-suspend the resulting spaces quaternionically.
This latter argument is delicate and uses nontrivial results from
the theory of cycles on quasi-projective varieties.  This is all
done in Section \ref{sec:6}.

\section{Stabilization and the ring structure}
\label{sec:3}

In this section we examine the limit over $V\subset \bbh^{\infty}$
of the cycles spaces $\cz_{\bbh}^q(V)$.  There are two series: $q$
even and $q$ odd, with different limits.  The algebraic join
induces a product on the homotopy groups of these spaces, and we
shall compute structure of the resulting $\bbz_2$--ring. We begin
with the following.

\begin{proposition}
\label{prop:3.1}
Let $V$ and $W$ be quaternionic vector spaces of complex
dimensions $2v$ and $2w$ respectively. Then for each $q$ with $0<
q < 2v$,  the inclusion
$$
\cz_{\bbh}^q(\bbp(V)) \subset \cz_{\bbh}^{q+2w}(\bbp(V\oplus
W))
$$
induces an injection on homotopy groups.
\end{proposition}
\begin{proof}
There are two cases to consider: $q$ even and $q$ odd. The
arguments are analogous.  In both cases  one reduces to 0--cycles
by quaternionic suspension: Corollary \ref{cor:6.4} in the even
case and Theorem \ref{thm:2.2} in the odd case. Then by applying
the Dold--Thom Theorem \cite{QDT} to 0--cycles, it suffices to show
that the maps in homology
$$ H_*(\bbp(V)/\bbz_2;\, \bbz) \arr H_*(\bbp(V\oplus W)/\bbz_2; \bbz)$$
and
$$H_*(Q(V)/\bbz_2;\, \bbz) \arr H_*(Q(V\oplus W)/\bbz_2; \bbz),$$
where $Q(V)$ is defined in \eqref{eq:6.2}, are injective.  This is
straightforward.
\end{proof}

\begin{corollary}
\label{cor:3.2}
 Consider the limiting spaces
$$\cz_{\bbh}^{\ev} = \lim_{n,q
  \rightarrow\infty}\cz_{\bbh}^{2q}(\bbp_{\bbc}(\bbh^n)) \aand
  \cz_{\bbh}^{\odd} = \lim_{n,q
  \rightarrow\infty}\cz_{\bbh}^{2q+1}(\bbp_{\bbc}(\bbh^n)).$$
There are canonical homotopy equivalences:
$$\aligned
\cz_{\bbh}^{\ev} &\cong \prod_{j=0}^{\infty} K(\bbz, 4j) \times
\prod_{j=0}^{\infty} K(\bbz_2, 4j+2)\\
\cz_{\bbh}^{\odd}  &\cong \prod_{j=0}^{\infty} K(\bbz, 4j)
\times \prod_{j=0}^{\infty} K(\bbz_2, 4j+1)
\endaligned
$$
\end{corollary}
\begin{proof}
This corollary follows from Proposition \ref{prop:3.1}, Theorem \ref{thm:2.3}, and 
\cite[Theorem~A.5]{QQLLM}.
\end{proof}

We now observe that the algebraic join gives   biadditive maps
$$
\# \co  \cz_{\bbh}^q(\bbp(V))\wedge \cz_{\bbh}^{q'}(\bbp(V')) \arr
\cz_{\bbh}^{q+q'}(\bbp(V\oplus V'))
$$
for all quaternionic vector spaces $V,V'$ and for all $q,q'$.
These maps induce pairings
$$
\pi_k\cz_{\bbh}^q(\bbp(V)) \otimes
\pi_{k'}\cz_{\bbh}^{q'}(\bbp(V')) \arr
\pi_{k+k'}\cz_{\bbh}^{q+q'}(\bbp(V\oplus V'))
$$
which together with Proposition \ref{prop:3.1} and Corollary
\ref{cor:3.2} above gives us the following.

\begin{proposition}
\label{prop:3.3}
Set $\grpi = \grpi_*^0 \oplus \grpi_*^1$, where
\begin{equation}
\label{eq:3.1}
\grpi_*^0 = \pi_*\cz_{\bbh}^{\ev} \aand \grpi_*^1 = \pi_*\cz_{\bbh}^{\odd}.
\end{equation}
Then the algebraic join gives $\grpi$ the structure of a
$\bbz_2{\times}\bbz$--graded ring.
\end{proposition}

The following determination  of this ring will be established in
Section~\ref{sec:7}.

\begin{theorem}
\label{thm:3.4}
The subring $\grpi_*^0$ admits an isomorphism
\begin{equation}
\label{eq:3.2}
\grpi_*^0 \cong  \bbz[x,u]/(2u, u^2)
\end{equation}
where $x$ corresponds to the generator of $\pi_4 \cz_{\bbh}^{\ev}
\cong \bbz$, $u$ corresponds to the generator of $\pi_2
\cz_{\bbh}^{\ev} \cong \bbz_2$, and $(2u, u^2)$ denotes the ideal
in the ring $\bbz[x,u]$ generated by $2u$ and $u^2$.

With respect to the isomorphism  \eqref{eq:3.2}, one has that
$u\cdot\grpi_*^1 = 0$ and $\grpi_*^1$ is the $\bbz[x]$--module
\begin{equation}
\label{eq:3.3}
\grpi_*^1  \cong \bbz[x] \lambda \oplus \bbz_2 [x] v
\end{equation}
where $\lambda$ corresponds to the generator of $\pi_0\grpi_*^1 =
\bbz$ and $v$ corresponds to the generator of $\pi_1\grpi_*^1 =
\bbz_2$.

The elements $\lambda$ and $v$ satisfy the relations
\begin{equation}
\label{eq:3.4}
\lambda^2 = 4, \qquad \lambda\cdot v = 0 \aand  v^2 = 0.
\end{equation}
\end{theorem}

\begin{note}
\label{note:3.5}
Note that $\grpi_*$ is a $\bbz_2$--graded $\bbz[x]$--algebra with
even generators $1$ and $u$ and ``companion'' odd generators
$\lambda$ and $v$.
\end{note}

\section{Extending functors from $K$--theory}
\label{sec:4}

As in the case of real cycles, certain basic  functors from
representation theory extend to quaternionic  algebraic cycles.
The constructions parallel those of the real case but curiously
the roles are interchanged.
\medskip

\noindent\textbf{Quaternionification}\qua To any complex vector space
$V$ we can associate the quaternionic vector space
$(V\otimes_{\bbc}\bbh, \jj)$ where
$$
(V\otimes_{\bbc}\bbh, \jj)  \equdef  V \oplus \overline{V}
\quad \text{and} \quad \jj(v,w) = (-w,v).
$$
Recall that if $V$ is a complex vector space, we define $\overline
V$ to be the same underlying real vector space with complex
structure changed from $J$ to $-J$. For any $q<\dim(V)$ we have a
map
\begin{equation}
\label{eq:4.1}
\cz^q_{\bbc}(\bbp(V)) \arr \cz^{2q}_{\bbh}(\bbp(V\oplus
\overline{V}))
\end{equation}
$$
c \mapsto c \# c .\leqno{\hbox{defined by}}
$$
This construction gives rise to commutative diagrams
\begin{equation}
\label{eq:4.2}
\bfig
\square<1000,400>[G^q_{\bbc}(\bbp(V))`G^{2q}_{\bbh}(\bbp(V\otimes_{\bbc}\bbh))`
  \cz^q_{\bbc}(\bbp(V))`\cz^{2q}_{\bbh}(\bbp(V\otimes_{\bbc}\bbh));
  `c`c_{\bbh}`]
\efig
\end{equation}
where $G^{2q}_{\bbh}$ denotes the Grassmannian of quaternionic
linear subspaces  of  quaternionic codimension $q$.  Diagram
\eqref{eq:4.2} stabilizes to a commutative diagram
\begin{equation}
\label{eq:4.3}
\bfig
\square<700,400>[BU_q`BSp_{q}`\cz^q_{\bbc}`\cz^{2q}_{\bbh};
  \quat`c`c_{\bbh}`\quatt]
\efig
\end{equation}
which in light of \cite[Theorem~1]{QL}, \cite[Theorem~A.5]{QQLLM} and
Theorem \ref{thm:2.3} can be rewritten canonically as
\begin{equation}
\label{eq:4.4}
\bfig
\square<1500,400>[BU_q`BSp_{2q}`\prod_{k=0}^q K(\bbz,
  2k)`\prod_{k=0}^{q} K(\bbz, 4k) \times \prod_{k=1}^{q} K(\bbz_2, 4k{-}2).;
  \quat`c`c_{\bbh}`\quatt]
\efig
\end{equation}
\medskip

\noindent\textbf{The forgetful functor}\qua Consider a quaternionic
vector space $(V,\jj)$ and the functor $(V,\jj) \mapsto V$ which
forgets the quaternionic structure. For any $q<\dim_{\bbh}(V)$ we
have a map
\begin{equation}
\label{eq:4.5}
 \cz^{2q}_{\bbh}(\bbp(V)) \arr \cz^{2q}_{\bbc}(\bbp(V))
\end{equation}
which simply includes the $\jj$--fixed cycles into the group of all
cycles. This gives commutative diagrams
\begin{equation}
\label{eq:4.6}
\bfig
\square<1000,400>[G^{2q}_{\bbh}`G^{2q}_{\bbc}(\bbp(V))`
  \cz^{2q}_{\bbh}(\bbp(V))`\cz^{2q}_{\bbc}(\bbp(V)).;
  `c_{\bbh}`c`]
\efig
\end{equation}
Note that $G^{2q}_{\bbh}$ is exactly the subset of $\jj$--fixed
planes $G^{2q}_{\bbc}$.  The diagram \eqref{eq:4.6}  stabilizes to
\begin{equation}
\label{eq:4.7}
\bfig
\square<700,400>[BSp_q`BU_{2q}`\cz^{2q}_{\bbh}`\cz^{2q}_{\bbc};
  \real`c_{\bbh}`c`\reall]
\efig
\end{equation}
where ${\real}\co BSp_q \to BU_{2q}$ is the map induced by the
standard embedding $Sp_q \to U_{2q}$ given by forgetting the
quaternionic structure.  Diagram \eqref{eq:4.7} can be rewritten
as
\begin{equation}
\label{eq:4.8}
\bfig
\square<1500,400>[BSp_q`BU_{2q}`\prod_{k=0}^{q} K(\bbz, 4k)
  \times\prod_{k=1}^{q} K(\bbz_2, 4k-2)`\prod_{k=0}^q K(\bbz, 2k);
  \real`c_{\bbh}`c`\reall]
\efig
\end{equation}
where these splittings are canonical; cf Theorem \ref{thm:2.3},
\cite[Theorem~1]{QL} and \cite[Theorem~A.5]{QQLLM}.

Let $\iota_{2k} \in H^{2k}(K(\bbz,2k);\,\bbz)$ be the fundamental
class as above, and denote by $\xi^q_{\bbh}$ the universal
quaternionic bundle of quaternionic rank $q$ over $BSp_q$.  Since
$c^*\iota_{2k}$ is the universal $k$th Chern class, the
commutativity of \eqref{eq:4.8} shows that
\begin{equation}
\label{eq:4.9}
 {c_{\bbh}}^*\Phi^*
(\iota_{2k}) = \phi^*c_k\left(\xi_{\bbc}^{2q}\right)  =
c_k\left(\xi_{\bbh}^{q}\right) = \begin{cases}
\sigma_m    &\text{if $k=2m$ and $0\leq m\leq q$},\\
0   &\text{otherwise,}
\end{cases}
\end{equation}
where $\sigma_1,\ldots,\sigma_q$ are the canonical generators of
$H^*(BSp_q;\, \bbz) = \bbz[\sigma_1,\ldots,\sigma_q]$.

The map $\Phi$ is entirely determined up to homotopy by the
following result.

\begin{theorem}
\label{thm:4.1}
Let   $\iota_{2k} \in H^{2k}(K(\bbz, 2k); \bbz) = \bbz$ be the
fundamental class.  Then
$$
\Phi^* \iota_{4k} = \iota_{4k}  \qquad\text{and}\qquad \Phi^*
{\iota}_{4k-2} = 0
$$
for all $k$.
\end{theorem}
\begin{proof}
The diagrams  \eqref{eq:4.8} successively embed one into the next
for $q=1,2, \ldots$ by taking linear embeddings $\dots\subset
\bbh^n \subset \bbh^{n+1}\subset \dots$.  Passing to quotients and
using Theorem \ref{thm:6.6}  below gives
\begin{equation}
\label{eq:4.10}
\bfig
\square<1500,400>[BSp_q`BU_{2q}`K(\bbz, 4q) \times  K(\bbz_2, 4q-2)`
  K(\bbz, 4q)\times  K(\bbz, 4q-2);
  \real`c_{\bbh}`c_q\times c_{q-1}`\widehat{\reall}]
\efig
\end{equation}
where ${\widehat{\reall}}\co \cz_{\bbh}^{2q}/\cz_{\bbh}^{2q-2}\arr
\cz_{\bbc}^{2q}/\cz_{\bbc}^{2q-2}$ is the induced map of quotient
groups. The second assertion of the theorem follows from the fact
that 
$$H^{4q-2}(K(\bbz, 4q) \times  K(\bbz_2, 4q-2) ;\,\bbz) = 0.$$
Now
$$
\aligned H^{4q}(K(\bbz, 4q) & \times  K(\bbz_2, 4q{-}2) ;\,\bbz) \\
&= H^{4q}(K(\bbz, 4q) ;\,\bbz) \oplus H^{4q}(K(\bbz_2, 4q{-}2) ;\,\bbz)\\
&= H^{4q}(K(\bbz, 4q) ;\,\bbz)  = \bbz \iota_{4q}.
\endaligned
$$
By \eqref{eq:4.9} we see that ${c_{\bbh}}^*\Phi^* (\iota_{2k})$ is
an additive generator and therefore $\Phi^* (\iota_{2k})$ must be
also.  This proves the first assertion.
\end{proof}
\medskip

\noindent\textbf{Relations}\qua  Consider the diagram
\begin{equation}
\label{eq:4.11}
\bfig
\square(0,0)|alla|/->`->`->`->/<1300,300>[\cz_{\bbh}^{q}`\cz_{\bbh}^{2q}`
  \prod_{j=0}^q K(\bbz,2j)`\prod_{j=0}^q K(\bbz,4j) \times
  \prod_{j=1}^{q} K(\bbz_2,4j{-}2);\Gamma`\cong`\cong`\Gamma]
\square(1300,300)|alra|/->`->`->`->/<1300,300>[BSp_{q}`BU_{2q}`
  \phantom{\cz_{\bbh}^{q}}`\cz_{\bbc}^{2q};``c`\Phi]
\morphism(0,600)/->/<0,-300>[BU_q`\phantom{\cz_{\bbh}^{q}};c]
\morphism(0,600)/->/<1300,0>[\phantom{BU_q}`\phantom{BSp_{q}};]
\morphism(1300,0)/->/<1300,0>[\phantom{\prod_{j=0}^q K(\bbz,4j) \times
  \prod_{j=1}^{q} K(\bbz_2,4j{-}2)}`\prod_{j=0}^{2q} K(\bbz,2j).;\Phi]
\morphism(2600,300)|r|/->/<0,-300>[\phantom{\cz_{\bbc}^{2q}}`
  \phantom{\prod_{j=0}^{2q} K(\bbz,2j).};\cong]
\efig
\end{equation}
From  \eqref{eq:4.1} we see  that if $V$ has a real structure
$\rho$, then under the  isomorphism $I\oplus \rho\co V \oplus
\overline V \arr V \oplus V$, the map  ${\Gamma} \co \cz_{\bbc}^{q}
\to \cz_{\bbh}^{2q}$ becomes ${\Gamma}(c) = c \# \rho_*(c)$.  It
follows that
$$
\Phi \circ{\Gamma}(c) = c \# \rho_*(c)
$$
for $c \in \cz_{\bbc}^{q}$. As in \cite[Proposition~5.1]{QQLLM} we
conclude the following.

\begin{proposition}
\label{prop:4.2}
Let
 $\iota_{2k} \in H^{2k}(K(\bbz, 2k); \bbz)
= \bbz$ be the fundamental class.  Then for each $k$ the
composition $\Phi \circ{\Gamma}$ satisfies
\begin{equation}
\label{eq:4.12}
(\Phi \circ{\Gamma})^* \iota_{2k} = \sum_{i+j=k}
(-1)^j\iota_{2i} \cup \iota_{2j}.
\end{equation}
\end{proposition}
\medskip

Note that in particular $(\Phi \circ{\Gamma})^* \iota_{4k-2}=0$ as
predicted by Theorem \ref{thm:4.1}.

Combining Theorem \ref{thm:4.1} and Proposition \ref{prop:4.2}
gives the following.

\begin{corollary}
\label{cor:4.3}
For each $k\leq q$ one has
$$
\G^*\iota_{4k} = \underset{i,j \leq q}{\sum_{i+j=2k}}
(-1)^j\iota_{2i}\cup \iota_{2j} + (-1)^{k} \iota_{2k}^{\, 2}.
$$
\end{corollary}
To completely determine   $\G$ up to  homotopy we need to compute
the classes $\G^* \widetilde{\iota}_{4k-2}$  where $
\widetilde{\iota}_{4k-2} \in H^{4k-2}(K(\bbz_2, 4k-2);\, \bbz_2) $
denotes the fundamental class. From the commutative diagram
\eqref{eq:4.4} we see that
$$
c^*\G^*\widetilde{\iota}_{4k-2} =
\g^*c_{\bbh}^*\widetilde{\iota}_{4k-2} = 0
$$
since $H^{4k-2}(BSp_q;\bbz_2) = 0$. Thus
$\G^*\widetilde{\iota}_{4k-2}$ lies in the kernel of $c^*$ on mod
2 cohomology. (Note that $c^*$ is injective on $\bbz_2[\iota_1,
\dots, \iota_q]$.)\ However, a complete calculation of this class
remains to be done.


\section{Infinite loop space structures.}
\label{sec:5}

In this section we carry the discussion in Section~6 of \cite{QQLLM}
over to the quaternionic case. Given a  quaternionic vector space
$(V,\jj)$ with $\dim_{\bbc}(V) = 2q$, we define  $\ci_*$--functors
$$
T_{G_{\bbh}}(V) = G_{\bbh}^{2q}(\bbp(V\oplus V)  \aand
T_{{Z}_{\bbh}}(V) = {\cz}_{\bbh}^{2q}(\bbp(V\oplus V).
$$
The  action on morphisms and the natural transformations
$\omega_{G_{\bbh}}$ and $\omega_{{Z}_{\bbh}}$ are defined exactly
as in \cite[Section~6]{QQLLM}. The inclusion
$$
T_{G_{\bbh}}(V) \subset T_{{Z}_{\bbh}}(V)
$$
as cycles of degree one is a natural transformation of
$\ci_*$--functors. As seen in \cite[page 16]{QMay}, the limiting
space
$$
\lim_{q\to\infty} T_{G_{\bbh}}(\bbh^q) = \mathbf{BSp}
$$
is a connected $\call$--space whose associated infinite loop space
structure coincides with the usual connective Bott structure. The
arguments of \cite[Section~6]{QQLLM} and \cite[Theorems~6.9 and
6.10]{QQLLM} apply directly to prove the following.

\begin{theorem}
\label{thm:5.1}
The limiting space $\cz^{\ev}_{\bbh}$ is an $E_{\infty}$--ring
space  and forms the 0--level  space  of an $E_{\infty}$--ring
spectrum.  The component $\cz^{\ev}_{\bbh}(1)$ consisting of
cycles of degree 1 carries an infinite loop space structure which
enhances the algebraic join and for which the induced mapping
$$
\mathbf{BSp}  \arr \cz^{ev}_{\bbh}(1)
$$
is a map of infinite loop spaces.
\end{theorem}

Consider now the  ``forgetful'' homomorphism
$$
\Phi\co \cz^{2q}_{\bbh}(\bbp(V\oplus V))\arr
\cz^{2q}_{\bbc}(\bbp(V\oplus V))
$$
defined in the last section.  This is a natural transformation of
$\ci_*$--functors, and so we have:

\begin{proposition}
\label{prop:5.2}
The limiting ``forgetful''  homomorphism
$$
\Phi\co \cz^{\ev}_{\bbh} \arr \cz^{\infty}_{\bbc}
$$
is a map of $E_{\infty}$ ring spaces. In particular, its
restriction $ \Phi\co \cz^{\ev}_{\bbh}(1)\arr
\cz^{\infty}_{\bbc}(1) $ is an infinite loop map.
\end{proposition}

Now the ``quaternionification'' maps
$$
\Gamma\co \cz^{2q}_{\bbc}(\bbp(V\oplus V)) \arr
\cz^{4q}_{\bbh}(\bbp(V\oplus \bar V\oplus V\oplus \bar  V))
$$
are {\sl not} additive mappings.  Nevertheless,  they do give
natural transformations of $\ci_*$--functors (with values in
topological spaces, not topological groups).  Hence we have:

\begin{proposition}
\label{prop:5.3}
The limiting ``quaternionification'' mapping
$$
\Gamma\co \cz^{\infty}_{\bbc} \arr \cz^{\ev}_{\bbh}
$$
is a  mapping of $\call$--spaces.  In particular, its restriction $
\Gamma\co \cz^{\infty}_{\bbc}(1) \arr \cz^{\ev}_{\bbh}(1) $ is an
infinite loop map.
\end{proposition}

\noindent These maps sit in a commutative diagram:
$$
\bfig
\square(0,0)<700,400>[\mathbf{BSp}`\mathbf{BU}`
  \cz^{\ev}_{\bbh}`\cz^{\infty}_{\bbc};\phi```\Phi]
\square(700,0)<700,400>[\phantom{\mathbf{BU}}`\mathbf{BSp}`
  \phantom{\cz^{\ev}_{\bbh}}`\cz^{\ev}_{\bbh};\gamma```\Gamma]
\efig
$$

\section{Proof of Theorem \ref{thm:2.3}}
\label{sec:6}

Part (ii) of this result was established  in \cite[Theorem~6.4]{QLLM},
so it remains only to prove part (i). The Quaternionic
Suspension Theorem of \cite[Theorem~6.1]{QLLM} applies also to cycles
of even codimension, but since quaternionic suspension changes
cycle-dimension by 2, one is unable in this case to reduce to
0--cycles where the Dold--Thom Theorem can be used.  We shall solve
this problem by ``replacing'' $\bbp(V)$ with an even-dimensional
variety $Q(V)$.

To begin   consider the  Veronese embedding $ v\co  \bbp(V)
\hookrightarrow \bbp(\sm(V)) $ which converts $\bbp(V)$ to a real
subvariety under a real structure $r\co \bbp(\sm(V)) \to
\bbp(\sm(V))$ coming from a complex conjugation on $\sm(V)$.  (To
see this explicitly choose coordinates $ z_\alpha + w_\alpha j$, for
$\alpha =1,\ldots,n$,  on $V = \bbh^n$ and note that
$v(z_\alpha,w_\alpha) = (z_\alpha z_\b, w_\alpha w_\b, z_\alpha
w_\b)$.) \ The following is a direct consequence of \cite{QLam}.
See also \cite{QLLM}.

\begin{proposition}
\label{prop:6.1}
The  $\bbz_2$--equivariant complex suspension map
\begin{equation}
\label{eq:6.1}
\csus\co  \cz^q(\bbp(V))  \arr \cz^q(Q(V))
\end{equation}
where
\begin{equation}
\label{eq:6.2}
 Q(V) \equiv \csus\{v(\bbp(V))\}  =
\operatorname{Thom}\left\{\cco{\bbp(V)}(2)\right\}
\end{equation}
is a $\bbz_2$--homotopy equivalence.
\end{proposition}
The idea now is to compute the $\bbz_2$--homotopy type of the
spaces $\cz^q(Q(V))$, for $q$ even, by ``de-suspending'' to the
case of 0--cycles. To begin we   fix some notation. Let
$$
\U^{2n} \equiv Q(\bbh^n) - \{\infty\} \equiv \cco{\pch n}(2)
 \xrightarrow{p} \pch n
$$
denote the square of the complex hyperplane bundle over $\pch n$.
Our real structure, $r\co Q(\bbh^n) \arr Q(\bbh^n)$, when restricted
to $\U^{2n}$, is an anti-linear bundle map for which the diagram
$$
\bfig
\square<700,400>[\U^{2n}`\U^{2n}`\pch n`\pch n;r`p`p`q]
\efig
$$
commutes.  (This bundle map is the one naturally induced on
$\cco{}(2)$ via multiplication by the quaternion $j$.)  Now the
topological groups of algebraic cycles, $\cz^q(U)$ are defined for
any quasi-projective variety $U$ (cf \cite[Definition~4.5]{QLi}), and,
since $Q(\bbh^n) - \U^{2n}$ consists of a single point, there are
$\bbz_2$--equivariant homeomorphisms
\begin{equation}
\label{eq:6.3}
 \cz^q(\U^{2n}) \xrightarrow{\cong}
\begin{cases}
 \cz^q(Q(\bbh^n)) & \ \text{for $q<2n$}, \\
 \cz^q(Q(\bbh^n))/\bbz &\  \text{for $q=2n$.}
\end{cases}
\end{equation}
We now observe that there is a commutative diagram of
$\bbz_2$--equivariant bundle maps
\begin{equation}
\label{eq:6.4}
\bfig
\square<1300,400>[\U^{2n+2}{-}\U^{2}`\bbp_{\bbc}(\bbh^n\oplus\bbh){-}\pch{}`
  \U^{2n}`\pch n;p`\pt`\pi`p]
\efig
\end{equation}
where $\pi$ is linear projection and $\pt$ is defined as follows.
Let $\lo \arr \pch n$ and $\ell \arr  \bbp_{\bbc}(\bbh^n\oplus
\bbh) $ denote the tautological complex line bundles $\cco{}(-1)$,
and note that
$$
\U^{2n+2} = \Hom(\ell\otimes\ell, \bbc) \qquad\text{and}\qquad
\U^{2n} = \Hom(\lo\otimes\lo, \bbc).
$$
The linear projection $\bbh^n\oplus \bbh\arr \bbh^n$ induces a
bundle mapping $\pi_* \co  \ell\arr\lo$ covering $\pi$ which is an
isomorphism on fibres.   The map $\pt$ in \eqref{eq:6.4} is given
by $\pt(h) = h \circ (\pi_*^{-1}\otimes \pi_*^{-1})$ for $h \in
\Hom(\ell\otimes\ell, \bbc)$.

Our main assertion here is the following.

\begin{proposition}
\label{prop:6.2}
The flat pull-back of cycles gives a \zt-homotopy equivalence
$$
\pt^* \co  \cz^q(\U^{2n}) \arr \cz^q(\U^{2n+2}-\U^{2})
$$
for all $q\leq 2n$.
\end{proposition}

\begin{inote}
\label{inote:6.3}
A quick proof of Proposition \ref{prop:6.2} can be given for
$q{<}2n$  as follows.  By \cite{QLLM} and \cite{QLam} the flat
pull-back of cycles gives equivariant homotopy equivalences
\begin{eqnarray*}
\pi^*\co \cz^q(\pch n) &\arr& \cz^q(\bbp_{\bbc}(\bbh^n\oplus \bbh)-\pch{})
\cong \cz^q(\pch {n+1}) \\
\text{and}\qquad p^*\co \cz^q(\pch {n}) &\arr& \cz^q(\U^{2n})
\end{eqnarray*}
for all $q$ and $n$; see also \cite[Proposition~4.15]{dSii}. Equivariant
excision arguments (cf \cite{QLi,QQLi,QLLM} or
\cite[Remark~4.14]{dSii}) then show that
$$
p^*\co \cz^q(\bbp_{\bbc}(\bbh^n\oplus \bbh)-\pch { })    \arr
\cz^q(\U^{2n+2}-\U^2)
$$
is also a \zt homotopy equivalence, and the Proposition follows
from the commutativity of \eqref{eq:6.4}.  Unfortunately this will
not allow us to reduce to 0--cycles, so we must construct a proof
directly in this case.
\end{inote}

\begin{proof}[Proof of Proposition~\ref{prop:6.2}]
We compactify each $\U^{2n} \cong \cco{\pch n}(2)$ by taking the
projective closure
$$
\ob^{2n}  = \bbp\{\cco{}(2)\oplus\bbc\}.
$$
This gives us a fibre square
\begin{equation}
\label{eq:6.5}
\bfig
\square<1200,400>[\ob^{2n+2}{-}\ob^{n}`\bbp_{\bbc}(\bbh^n\oplus \bbh){-}\pch{}`
  \ob^{2n}`\pch n;p`\pt`\pi`p]
\efig
\end{equation}
of smooth \zt maps.  We recover the diagram \eqref{eq:6.4} by
removing the restriction of $\pt$ to the ``$\infty$--section''
$\bbp_{\infty}^{2n-1} \subset \ob^{2n}$. Observe that taking the
graph gives an equivariant isomorphism
\begin{equation}
\label{eq:6.6}
\bfig
\square/->`->`->`=/<1000,400>[\Hom_{\bbc}(\lo,\bbh)`
  \bbp_{\bbc}(\bbh^n\oplus\bbh){-}\pch{}`\pch n`\pch n;
  \cong`\pi`\pi`]
\efig
\end{equation}
where the  \zt action on $\Hom_{\bbc}(\lo, \bbh)$ is given by
sending a linear map $h\co \lo \arr \bbh$ to $j(h) = j\circ h \circ
j^{-1}$.  Via this isomorphism we can rewrite the fibre square in
\eqref{eq:6.5} as a pull-back diagram
\begin{equation}
\label{eq:6.7}
\bfig
\square<1000,400>[\Hom_{\bbc}(\plo, \bbh)`\Hom_{\bbc}(\lo, \bbh)`
  \ob^{2n}`\pch n.;p`\pt`\pi`p]
\efig
\end{equation}
The restriction
$$
\bfig
\square/@^{(->}`->`->`@^{(->}/<1000,400>[\Hom_{\bbc}(\plo,\bbh)`\Hom_{\bbc}(\plo,\bbh)`
  \bbp_{\infty}^{2n+1}`\ob^{2n};`\pt`\pt`]
\efig
$$
of $\pt$ to the infinity-section is isomorphic (via $p$) to the
bundle
$$\bbp_{\bbc}(\bbh^n\oplus \bbh)-\pch { } \arr \pch n$$
for
which the Quaternionic Suspension Theorem holds.  Thus by
equivariant excision (see \cite[Remark~4.14]{dSii}) our
Proposition will follow if we can prove that
$$
\pt^* \co  \cz^q(\ob^{2n}) \arr \cz^q(\Hom_{\bbc}(\plo, \bbh))
$$
is a \zt homotopy equivalence for all $q$.  For our application
the interesting case (and by Note \ref{inote:6.3} the only
remaining case) is where $q = 2n$. Thus we shall prove the
assertion that
\begin{equation}
\label{eq:6.8}
 \pt^* \co  \cz_0(\ob^{2n}) \arr \cz_2(\Hom_{\bbc}(\plo, \bbh))
\ \text{is a \zt homotopy equivalence}
\end{equation}
where $\cz_p$ denotes the group of cycles of {\sl dimension}
$p$.

To prove \eqref{eq:6.8} we consider the submonoid
\begin{equation}
\label{eq:6.9}
\ct^+_2 \subset \cc_2(\Hom_{\bbc}(\plo, \bbh)),
\end{equation}
of effective 2--cycles which meet the zero-section in proper
dimension (namely 0), and denote by
\begin{equation}
\label{eq:6.10}
 \ct_2 \subset \cz_2(\Hom_{\bbc}(\plo, \bbh))
\end{equation}
the induced homomorphism of \naive topological group completions.

Observe now that scalar multiplication
$$
\Phi_t \co  \Hom_{\BC}(\plo, \BH) \arr\Hom_{\BC}(\plo, \BH)
$$
by real numbers $t>0$   gives bundle maps commuting with the
$\BZ_2$--action, and pulling to the normal cone gives a
$\BZ_2$--deformation retraction
$$
\ct_2 \arr \pt^*\left\{\cz_0({\ob^n})\right\}
$$
(cf \cite[Assertion 1, Section~2.3]{QLLM}, \cite{QFL}, \cite{QL}).
Therefore, it remains only to show that the inclusion
\eqref{eq:6.10} is a \zt homotopy equivalence.

We shall proceed in analogy with the arguments in
\cite[Section~2.3]{QLLM}. We consider the direct sum
$$
\Hom_{\BC}(\plo, \BH)\oplus \Hom_{\BC}(\plo, \BH) \arr  {\ob^n}
$$
and   choose two distinct projections
$$
\pi_0, \pi_{\infty} \co  \Hom_{\BC}(\plo, \BH)\oplus \Hom_{\BC}(\plo,
\BH) \arr \Hom_{\BC}(\plo, \BH).
$$
The map    $\pi_0$ is simply projection onto the first factor, and
$\pi_{\infty} = \pi_0 \circ S$ where
$$S = \begin{pmatrix}
\Id & \epsilon J \\
0 & \Id
\end{pmatrix},$$
$\epsilon >0$, and $J\co  \{0\}\oplus\Hom_{\BC}(\plo, \BH)
\xrightarrow{\cong} \Hom_{\BC}(\plo, \BH)\oplus \{0\} $ is the
canonical isomorphism between the second and first factors.

These projections can be viewed alternatively as follows.  Via
\eqref{eq:6.6} and \eqref{eq:6.7} we obtain a pull-back diagram
$$
\bfig
\square<1500,400>[\Hom_{\BC}(\plo,\BH)\oplus\Hom_{\BC}(\plo,\BH)`
  \BP_{\BC}(\BH^{n+2})-\BP_{\BC}(\BH^2)`
  \ob^{2n}`\BP_{\BC}(\BH^n),;p`\pt\oplus\pt`\pi_0\oplus\pi_\infty`p]
\efig
$$
and  $\pi_0, \pi_\infty$ are just the pull-backs of the
 ones constructed  in $\BP_{\BC}(\BH^{n+2})$ by projecting away
from   quaternionic lines $\lambda_0, \lambda_{\infty} \subset
\bbp_{\bbc}(\bbh^2)$. \medskip

Consider now the  open dense  set $\cu(d)$ of divisors $D$ of
degree $d$ on $\BP_{\BC}(\BH^{n+2})$ with the property that $D$
meets $jD$ in proper dimension and that  $\BD = D\bullet jD$ (and
all scalar multiples $t\BD$ for $0<t\leq 1$) do not meet the
vertices $\lambda_0, \lambda_{\infty}$ of our projections To each
$D \in \cu(d)$ we associate the pull-back cycle $\tbbd = p^*\bbd$
and define a transformation
$$
\Psi_D \co  \cc_2(\Hom_{\BC}(\plo, \BH)) \arr
\cc_2(\Hom_{\BC}(\plo, \BH))
$$
as in \cite[(2.6.1)]{QLLM}, \cite[page 285]{QL} by setting
$$
\Psi_D(c)  = (\pi_{\infty})_*\left\{\pi_0^* c \bullet
\tbbd\right\}.
$$
Note that $\pi_{\infty}$ and $\pi_0$ are proper on $\tbbd$. Note
also that if $\deg(D) = d$, then $\lim_{t\to 0}\Psi_{tD} = d^2\cdot
\Id$.

The arguments given in \cite[pp 634--641]{QLLM} and the proof of
\cite[Theorem~6.1]{QLLM}  may now be repeated in this context. The
important point is to show that there is a function $N(d)$ such
that $N(d) \to \infty$ as $d\to \infty$ and with the property that
for any irreducible subvariety $Z\subset \Hom_{\BC}(\plo, \BH)$ of
dimension 4,
$$
\codim \left( \gB_Z(d)\right) \geq N(d)
$$
$$
\gB_Z(d) = \{D\in \cu(d) : \dim(Z\cap \tbbd ) \geq 3\}
\leqno{\hbox{where}}
$$
is the set of ``bad'' divisors of degree $d$ for $Z$.

Now given such a $Z$, consider the irreducible subvariety $p(Z)
\subset \BP_{\BC}(\BH^{n+2})$. Note that
$$
\text{dim} (p(Z)) = \text{either 3 or 4}
$$
and that the fibres of $p\co  Z \arr  p(Z)$ are of dimension at most
1 since the fibres of $p$ are complex lines.

Suppose that $\dim(p(Z)) = 3$. In \cite[Lemma~2.7]{QLLM} it is
proved that, since the $\bbz_2$--action is free, there is a
function $N(d)$ independent of $Z$ and going to infinity with $d$
such that
\begin{equation}
\label{eq:6.11}
\codim \{D\in \cu(d) : \dim(p(Z)\cap \bbd) \geq 2\} \geq
N(d).
\end{equation}
 Since the fibre-dimension of $p$ is $\leq 1$, we see that
$$
\dim(p(Z)\cap \bbd) = 1\quad \Rightarrow\quad \dim
(Z\cap \tbbd) \leq 2
$$
and so the set of divisors in \eqref{eq:6.11}  contains
$\gB_Z(d)$.

Suppose now that $\dim(p(Z)) = 4$.  Then the generic fibre of $p\co  Z
\arr  p(Z)$ has dimension 0, and there is a subvariety $\Sigma
\subset p(Z)$ of dimension $\leq 3$ where the fibre dimension is
1.  Again from \cite[Section~2.7]{QLLM} we know that there is a
function $N(d)$ as above such that
\begin{equation}
\label{eq:6.12}
\codim \{D\in \cu(d) : \dim(p(Z)\cap \bbd) \geq 3, \ \text{or}
\ \dim(\Sigma\cap \bbd) \geq 1\} \geq N(d).
\end{equation}
The set of divisors in \eqref{eq:6.12} contains $\gB_Z(d)$, and so
we have proved the desired estimate on the codimension of the bad
sets.  The arguments of \cite[pp 640--641]{QLLM} carry through to
establish assertion \eqref{eq:6.8}, thereby completing the proof
of Proposition \ref{prop:6.2}.
\end{proof}

\begin{corollary}
\label{cor:6.4}
There are \zt homotopy equivalences
$$
\cz^q(Q(V)) \cong \cz^q(Q(V\oplus \bbh))
$$
for all $q \leq \dim_{\bbc}(V)$.
\end{corollary}
\proof
Note that if $q<2n$, then by \eqref{eq:6.3}
$$
\aligned \cz^{q}(Q(\bbh^n)) &= \cz^{q}(\U^{2n}) \cong
\cz^{q}(\U^{2n+2}-\U^2)
\equdef \cz^{q}(\U^{2n+2})/\cz_{2n+2-q}(\U^2) \\
&= \cz^{q}(\U^{2n+2}) = \cz^{q}(Q(\bbh^{n+1})).
\endaligned
$$
This case follows also from \eqref{eq:6.1} and \cite{QLLM}.  For
the final case note that
$$\begin{aligned}
\cz^{2q}(Q(\bbh^q))/\bbz &= \cz^{2q}(\U^{2q}) \cong
\cz^{2q}(\U^{2q+2}-\U^2)
= \cz^{2q}(\U^{2q+2})/\cz_{2}(\U^2) \\
&= \cz^{2q}(\U^{2q+2})/\bbz  = \cz^{2q}(Q(\bbh^{q+1}))/\bbz.
\end{aligned}\eqno{\qed}
$$

From \eqref{eq:6.1} and Corollary \ref{cor:6.4} we conclude that
there is a \zt homotopy equivalence
$$
\cz^{2q}(\pch n) \cong \cz_0(Q(\bbh^q))
$$
and therefore there is a homotopy equivalence
\begin{equation}
\label{eq:6.13}
\cz^{2q}_{\bbh}(\pch n) \cong \cz_0(Q(\bbh^q))^{\text{fixed}}.
\end{equation}
Since $\bbz_2$ acts freely outside of one point on $Q(\bbh^q)$,
there is a group isomorphism
\begin{equation}
\label{eq:6.14}
\cz_0(Q(\bbh^q))^{\text{fixed}}/\cz_0(Q(\bbh^q))^{\text{av}} =
\bbz_2
\end{equation}
from which it follows that $\pi_*(\cz_0(Q(\bbh^q))^{\text{fixed}})
= \pi_*(\cz_0(Q(\bbh^q))^{\text{av}})$. We also have  that
$$
\cz_0(Q(\bbh^q))^{\text{av}} = \cz_0(Q(\bbh^q)/\bbz_2).
$$
Now by the work of Dold--Thom \cite{QDT} we know that for a
connected finite complex $A$ there is a homotopy equivalence
$$
\cz_0(A) \cong \prod_{j\geq 0} K(H_j(A;\bbz), j).
$$
Therefore the first part of Theorem \ref{thm:2.3}(i)  follows from
the next Proposition.

\begin{proposition}
\label{prop:6.5}
$$
H_j( Q(\bbh^q)/\bbz_2;\, \bbz) = \begin{cases}
\bbz    &  \qquad  \text{if $j\equiv 0$ (mod 4) and $j\leq 4q$}, \\
\bbz_2     &  \qquad  \text{if $j\equiv 2$ (mod 4) and $j\leq 4q-2$}, \\
0      &  \qquad  \text{otherwise}.
\end{cases}
$$
\end{proposition}
\begin{proof}
Set $Y = Q(\bbh^q)/\bbz_2$ and $X = \pch q / \bbz_2$, and note
that $Y$ is the Thom space of a nonorientable real 2--plane bundle
$L \to  X$. ($L$ is simply the quotient of $\U^{2q} =
\cco{\bbp^{2q-1}}(2)$ by the quaternion involution.)\ Thus,
$$
H_*(Y;\, \bbz) = H_*(B_L, S_L;\, \bbz)
$$
where $S_L \subset B_L$ denote the unit circle and unit disk
bundles of $L$. By looking at the Hopf fibration one can see that
$$
S_L = S^{4q-1}/\bbz_4
$$
where $\bbz_4$ is generated by multiplication by the quaternion
$j$ on the unit sphere $S^{4q-1}\subset \bbh^q$.  We know that:
\begin{align*}
H_j( X ;\, \bbz) &= \begin{cases}
\bbz    &  \qquad  \text{if $j\equiv 0$ (mod  4) and $j\leq 4q-4$}, \\
\bbz_2     &  \qquad  \text{if $j\equiv 1$ (mod  4) and $j\leq 4q-3$}, \\
0      &  \qquad  \text{otherwise}
\end{cases} \\
H_j( S_L ;\, \bbz) &= \begin{cases}
\bbz    &  \qquad  \text{if $j = 0$ or $4q-1$}, \\
\bbz_4     &  \qquad  \text{if $j$ is odd and $j< 4q-1$}, \\
0      &  \qquad  \text{otherwise}
\end{cases}
\end{align*}
The long exact sequence in homology for the pair $(B_L, S_L)$
gives $H_{4q}(Y;\, \bbz) = \bbz$, and for $i < 2q$ it gives exact
sequences
\begin{equation}
\label{eq:6.15}
 0 \rightarrow H_{2i}(X) \rightarrow H_{2i}(Y) \rightarrow \bbz_4
   \rightarrow H_{2i-1}(X) \rightarrow H_{2i-1}(Y) \rightarrow 0.
\end{equation}
When $i=2k$ we get
$$
0 \to \bbz  \to H_{4k}(Y) \to \bbz_4 \to 0 \to H_{4k-1}(Y) \to 0,
$$
and so  $H_{4k-1}(Y)  = 0$ and we have the short exact sequence
\begin{equation}
\label{eq:6.16}
 0  \arr  \bbz  \arr H_{4k}(Y)  \arr \bbz_4  \arr 0.
\end{equation}
To understand this extension we consider the $\bbz_2$--homology
groups
$$
H_j( X ;\, \bbz_2) = \begin{cases}
0    &  \qquad  \text{if $j\equiv 3$ (mod  4) and $j\leq 4q-2$}, \\
\bbz_2      &  \qquad  \text{otherwise}.
\end{cases}
$$
By the Thom isomorphism,
$$
H_{j+2}( Y ;\, \bbz_2) \cong H_{j}( X ;\, \bbz_2)
$$
for all $j$.     In particular we have $H^{4k}(Y;\, \bbz_2) =
H_{4k}(Y;\, \bbz_2) = \bbz_2$.  We conclude that
$$
H_{4k}(Y;\, \bbz) = \bbz,
$$
since all other possibilities for the extension \eqref{eq:6.16}
force the dimension of the vector space $H^{4k}(Y;\, \bbz_2)$ to
be greater than $1$.

Now when  $i = 2k+2$ in \eqref{eq:6.15}  we have
$$
0   \to H_{4k+2}(Y) \to \bbz_4 \xrightarrow{\sigma} \bbz_2 \to
H_{4k+1}(Y) \to 0,
$$
and it remains to show that $\sigma \neq 0$.  For this we must
also consider the $\bbz_2$--homology groups
$$
H_j( S_L ;\, \bbz_2) = \begin{cases}
0    &  \qquad  \text{if $j > 4q-1$}, \\
\bbz_2      &  \qquad  \text{otherwise}.
\end{cases}
$$
From the pair $(B_L, S_L)$ we get the   exact sequence of
$\bbz_2$--homology groups
\begin{multline*}
0 \arr H_{4k+3}(Y) \arr H_{4k+2}(S_L) \arr H_{4k+2}(X) \arr \\
  \arr H_{4k+2}(Y) \arr H_{4k+1}(S_L) \arr H_{4k+1}(X) \arr  0
\end{multline*}
which becomes
$$
0 \to \bbz_2 \xrightarrow{\cong}  \bbz_2 \xrightarrow{0} \bbz_2
\xrightarrow{\cong} \bbz_2 \xrightarrow{0} \bbz_2
\xrightarrow{\cong} \bbz_2 \to 0.
$$
The generator of $H_{4k+1}(S_L;\, \bbz) =\bbz_4$ goes to the
generator of $H_{4k+1}(S_L;\, \bbz_2) =\bbz_2$, and from the line
above we see that
$$
H_{4k+1}(S_L;\, \bbz_2) \xrightarrow{\sigma\otimes \bbz_2}
H_{4k+1}(X;\, \bbz_2)
$$
is not zero.  Hence, $\sigma$ is not zero, and the proposition is
proved.
\end{proof}

We now observe that all the constructions in the proofs of
Propositions \ref{prop:6.1} and \ref{prop:6.2} preserve the
subgroups of averaged cycles, and so the suspension maps induce
homotopy equivalences of these subgroups. Therefore, the
suspension maps induce homotopy equivalences of the quotients
$$
\cz^{2q}(Q(\bbh^n))^{\text{fixed}}/\cz^{2q}(Q(\bbh^n))^{\text{av}}
 \cong
\cz_0(Q(\bbh^q))^{\text{fixed}}/\cz_0(Q(\bbh^q))^{\text{av}},
$$
and by \eqref{eq:6.14} the right hand side is the space of two
points. This gives the second half of Theorem \ref{thm:2.3}(i).

The canonical nature of the homotopy equivalences in Theorem
\ref{thm:2.3} is established in \cite[Appendix~A]{QQLLM}. \qed

From the argument above we can deduce the following.

\begin{theorem}
\label{thm:6.6}
Fix  $q<n$. Let
\begin{equation}
\label{eq:6.17}
 \cz_{\bbh}^{2q-2}(\bbp_{\bbc}(\bbh^{n-1})) \subset
\cz_{\bbh}^{2q}(\bbp_{\bbc}(\bbh^{n}))
\end{equation}
be the subgroup of cycles contained in the linear subspace
$\bbp_{\bbc}(\bbh^{n-1})$, and let
\begin{equation}
\label{eq:6.18}
\cz_{\bbh}^{2q-1}(\bbp_{\bbc}(\bbh^{n-1})) \subset
\cz_{\bbh}^{2q+1}(\bbp_{\bbc}(\bbh^{n}))
\end{equation}
be defined similarly.  Then the inclusion \eqref{eq:6.17} is
$(4q-3)$--connected, and \eqref{eq:6.18} is $(4q-1)$--connected.
Furthermore, there are canonical homotopy equivalences
$$\aligned
\cz_{\bbh}^{2q}/\cz_{\bbh}^{2q-2} &\cong K(\bbz,4q)\times K(\bbz_2, 4q-2), \\
\cz_{\bbh}^{2q+1}/\cz_{\bbh}^{2q-1} &\cong K(\bbz,4q)\times
K(\bbz_2, 4q-1).
\endaligned$$
\end{theorem}
\begin{proof}
As in the proof of \cite[Proposition~8.1]{QQLLM}   we see that
the short exact sequence
$$
0 \arr \cz^{2q-1}_{\bbh} \arr \cz^{2q+1}_{\bbh} \arr
\cz^{2q+1}_{\bbh}/\cz^{2q-1}_{\bbh}  \arr 0
$$
is a fibration sequence.  Quaternionic algebraic suspension
\cite[Theorem~6.1]{QLLM} shows that this sequence is equivalent to
the fibration sequence
\begin{equation}
\label{eq:6.19}
 0 \arr \cz_0(\bbp_{\bbc}^{2q-1}/\bbz_2) \xrightarrow{i_*}
\cz_0(\bbp_{\bbc}^{2q+1}/\bbz_2) \arr
\cz_0((\bbp_{\bbc}^{2q+1}/\bbp_{\bbc}^{2q-1})/\bbz_2)  \arr 0.
\end{equation}
One sees directly that $i\co  \bbp_{\bbc}^{2q-1}/\bbz_2 \subset
\bbp_{\bbc}^{2q+1}/\bbz_2$ induces an isomorphism of $H_k(
\bullet; \bbz)$ for all $k\leq 4q-1$.  Hence by the Dold--Thom
Theorem $i_*$ induces an isomorphism of $\pi_k(\bullet)$ for $k$
in the same range.  This proves the theorem for
$\cz^{2q+1}_{\bbh}$. The result for $\cz^{2q}_{\bbh}$ is proved
analogously using the suspension arguments above.
\end{proof}


\section{Proof of Theorem \ref{thm:3.4}}
\label{sec:7}

We begin by observing that for quaternionic vector spaces $V,V'$
there is a commutative diagram
$$
\bfig
\square<1500,400>[\cz^{q}_{\bbh}(\bbp(V)) \wedge \cz^{q'}_{\bbh}(\bbp(V'))`
  \cz^{q+q'}_{\bbh}(\bbp(V\oplus  V'))`
  \cz^{q}_{\bbc}(\bbp(V)) \wedge  \cz^{q'}_{\bbc}(\bbp(V'))`
  \cz^{q+q'}_{\bbc}(\bbp(V\oplus  V'));
  \#`\Phi \wedge \Phi`\Phi`\#]
\efig
$$
where the vertical maps $\Phi$ are given by the simple inclusion
of the quaternionic cycles into the group of all cycles (cf
Section \ref{sec:4}). Under stabilization these maps yield a {\sl
ring homomorphism}
\begin{equation}
\label{eq:7.1}
\Phi_* \co  \grpi_*   \arr \pi_*\cz_{\bbc}^{\infty}  = \bbz[s].
\end{equation}

\begin{proposition}
\label{prop:7.1}
Let $p \co  \bbp_{\bbc}(\bbh^2) \to S^4 = \bbp_{\bbh}(\bbh^2)$ be the
``Hopf mapping'' which assigns to a complex line the quaternion
line containing it.  Define
$$
f\co S^4 \to \cz^2_{\bbh}\left(\bbp_{\bbc}(\bbh^2)\right)
$$
by setting $f(\ell) = p^{-1}(\ell) = $ ``$\ell$''.  Then $[f] = x
\in \pi_4\cz_{\bbh}^{\ev} \cong \bbz$ is the generator.
Furthermore, under the ring homomorphism \eqref{eq:7.1}  one has
$$
\Phi_*(x) = s^2
$$
and so $ \Phi_*(x^m) = s^{2m} $ for all $m$. Similarly, let
$\lambda$ denote the generator of $\pi_0\cz^{\odd}_{\bbh} \cong
\bbz$. Then $ \Phi_*(x^m\lambda) = 2s^{2m} $ for all $m$. In
particular, $x^m$ and $x^m\lambda$ are additive generators for all
$m$.
\end{proposition}
\begin{proof}
Under the composition
$$
\pi_4 \cz^2_{\bbh}(\bbp^3_{\bbc}) \xrightarrow{\Phi}\pi_4
\cz^2_{\bbc}(\bbp^3_{\bbc})
 \xrightarrow{\cong} H_6( \bbp^3_{\bbc};\, \bbz) = \bbz
$$
the class of $[f]$ goes to class of its ``trace'', which is the
union of the lines parameterized by $f$ (cf \cite[Section~9]{QQLLM}).
This trace is exactly $\bbp^3_{\bbc}$ whose class is the generator
of $H_6$.  It follows that $\Phi_*(x^m)$ must be the generator
$s^2$.

By Theorem \ref{thm:2.3}(ii) we know that every cycle in
$\cz_{\bbh}^{2q+1}$ is deformable to an averaged cycle, the
homomorphism $\deg\co  \pi_0\cz_{\bbh}^{\odd} \to \bbz$ given by
projective degree has
$$
\operatorname{Im}\{\deg\} = 2\bbz.
$$
Now fix $x_0 \in \bbp_{\bbc}(\bbh)$  and consider  the cycle $c_0
= x_0+\jj x_0$ whose component generates $\pi_0\cz_{\bbh}^{\odd}$.
Define $\widetilde{f} \co  S^4 \to \cz_{\bbh}^3(\bbp^5_{\bbc})$   by
$\widetilde{f}(t) = c_0\# f(t)$.  Since  $\widetilde{f}$ is
homotopic to twice the suspension of $f$ we see that
$\widetilde{f}$ has image $2\bbz$. The same holds for
$\widetilde{f}\wedge f \wedge \ldots \wedge f$.
\end{proof}

\begin{proposition}
\label{prop:7.2}
Let $1$ denote the generator of $\pi_0\cz_{\bbh}^{\ev} = \bbz$ and
let $u$ and $v$ denote the generators of $\pi_2\cz_{\bbh}^{\ev}  =
\bbz_2$ and $\pi_1\cz_{\bbh}^{\odd}  = \bbz_2$ respectively. Then
$1$ is the multiplicative unit and the following relations hold in
the ring $\grpi_*$:
\begin{equation}
\label{eq:7.2}
\aligned
u\cdot \lambda = 0,\quad v\cdot \lambda &= 0,\quad u^2 = 0,\\
v^2 = 0,\quad   u\cdot v &= 0,\quad  \lambda^2 = 4.
\endaligned
\end{equation}
\end{proposition}
\begin{proof}
That $1$ is a multiplicative unit is an immediate consequence of
the Quaternionic Suspension Theorem \cite[Theorem~6.1]{QLLM}. From
the paragraph above we see that $\lambda^2$ is represented by
$c_0\# c_0$ which has degree $4$.  Since $\deg\co  
\pi_0\cz_{\bbh}^{\ev}  \to \bbz$ is an isomorphism we conclude
that $\lambda^2 = 4$. That $v^2 = 0$ is Lemma \ref{lem:7.3} below.
All the remaining relations are trivial.
\end{proof}

It remains to prove that $x^mu$ and $x^mv$ are additive generators
for all $m>0$.  To do this we will need explicit representatives
for these classes.

\def\jl{\ell^{\perp}}
Recall the isomorphisms
$$\pi_1\cz_{\bbh}^{\odd} =
\pi_1\cz_{\bbh}^{1}(\bbp_{\bbc}(\bbh)) =
\pi_1\cz_0(\bbp_{\bbc}/\bbz_2) = H_1(\bbp_{\bbc}/\bbz_2;\bbz)   =
H_1(\bbp_{\bbr};\bbz) =\bbz_2.$$
Unravelling these isomorphisms
one sees that the generator $v$ of $\pi_1\cz_{\bbh}^{\odd}$ is
represented by the following map.  Let $\ell\co [0,\pi] \to S^2$ be
the standard longitudinal curve joining the north and south poles.
Under the identification $S^2 = \bbp_{\bbc}^1$, we consider
$\ell(t)$ to be a complex line in $\bbc^2 = \bbh$. Note that for
$\ell \in \bbp_{\bbc}^1 \subset \cz_0(\bbp_{\bbc}^1) =
\cz_{\bbh}^{1}(\bbp_{\bbc}(\bbh))$ we have that
\begin{equation}
\label{eq:7.3}
\jj(\ell) =  \jl = \ \text{the antipodal image of the point}\ \ell.
\end{equation}
 We now define
$ \phi\co S^1 \to \pi_1\cz_{\bbh}^{1}(\bbp_{\bbc}(\bbh)) $ by
\begin{equation}
\label{eq:7.4}
 \phi(t) = \ell(t) + \jl(t) - (\ell_0+ \jl_0).
\end{equation}
where $\ell_0 = \ell(0)$. By \eqref{eq:7.3} we see that $\phi$ is
a map into $\jj$--averaged cycles and $\phi(0) = \phi(\pi)$.

\begin{lemma}
\label{lem:7.3}
The map $\phi$  represents the generator $v$ of
$\pi_1\cz_{\bbh}^{\odd}$, and $v^2=0$  in
$\pi_2\cz_{\bbh}^{\odd}$.
\end{lemma}
\begin{proof}
The first statement follows from the paragraph above. For the
second note that
$$\aligned
\phi(s)\#\phi(t) &= \bigl\{\ell(s) +\jl(s)  - (\ell_0+
\jl_0)\bigr\} \#
\bigl\{\ell(t) +\jl(t)  - (\ell_0+ \jl_0)\bigr\}\\
&= \bigl\{(\ell(s) - \ell_0)\#(\ell(t) - \ell_0) +
(\jl(s) - \jl_0)\#(\jl(t) - \jl_0)\bigr\} +\\
&\qquad\qquad \bigl\{(\ell(s) - \ell_0)\#(\jl(t) - \jl_0) +
(\jl(s) - \jl_0)\#(\ell(t)
- \ell_0)\bigr\}\\
\endaligned$$
$$\aligned
&= \bigl\{(\ell(s) - \ell_0)\#(\ell(t) - \ell_0) +
\jj(\ell(s) - \ell_0)\#(\ell(t) - \ell_0)\bigr\} +\\
&\qquad\qquad \bigl\{(\ell(s) - \ell_0)\#(\jl(t) - \jl_0) +
\jj(\ell(s) - \ell_0)\#(\jl(t) - \jl_0)\bigr\}  \\
&\equdef A(s,t) + B(s,t).
\endaligned$$
It is straightforward to see that $A$ and $B$ are homotopic as
maps from $S^2$ into $\cz_{\bbh}^2(\bbp_{\bbc}(\bbh^2))$.  Hence,
$v^2 = 2[A] = 0$ in $\pi_2\cz_{\bbh}^{\ev} = \bbz_2$ as claimed.
\end{proof}

\begin{proposition}
\label{prop:7.4}
Let
$$
F = \phi\# f \# \dots \# f \co  S^1\wedge  S^4\wedge \dots \wedge S^4
\arr \cz_{\bbh}^{2m+1}(\bbp_{\bbc}(\bbh^{2m+1})),
$$
where $f$ and $\phi$ are the maps from Proposition \ref{prop:7.1}
and \eqref{eq:7.4} respectively.  Then $[F] = v\cdot x^m$ is the
generator of $\pi_{4m+1}\cz_{\bbh}^{\odd} = \bbz_2$.
\end{proposition}

\begin{proof} \def\tF{\widetilde F} It suffices to show that $[F] \neq 0$. For
this we will consider the graph of $F$ as in \cite[9.12ff]{QQLLM}. To begin note that $F$  has the form
\begin{multline*}
F(t,\la_1,\dots,\la_m) \\ = (\ell(t)+\jl(t)-\ell_0-\jl_0)
\# (\la_1-\la_0) \# (\la_2-\la_0) \#\dots \# (\la_m-\la_0) \\
= \tF(t, \la_1,\dots,\la_m) + \jj \tF(t, \la_1,\dots,\la_m)
\end{multline*}
where $\tF = (\ell(t)-\ell_0) \# (\la_1-\la_0) \# \dots \#
(\la_m-\la_0)$ can be considered as a map into cycles on $X
\equdef \bbp_{\bbc}(\bbh^{2m+1})/\bbz_2$. Now by Theorem
\ref{thm:2.3}(ii) if $F$ is homotopic to 0, then it is homotopic
to zero through $\jj$--averaged cycles, and therefore $\tF$ is
homotopic to 0 as a map into cycles on $X$.

To show $\tF$ is not homotopic to 0 we consider its graph
 $\Gamma$  in  $S^{4m+1}\times X$ and
show that its projection $[\text{pr}_*\Gamma] \neq 0$ in
$H_{8m+1}(X;\, \bbz_2)$ (see \cite[Lemma~9.12]{QQLLM}). Now we
see that
$$
\text{pr}_*\Gamma =  \bigcup_{t, \la_1,\dots,\la_m}\ell(t) \#
\la_1\# \dots \# \la_m + \epsilon \equiv G + \epsilon
$$
where $\epsilon$ consists of terms which have dimension strictly
less than $8m+1$ and can be ignored.

Suppose now that $q = (q_0,q_1,\ldots,q_m) \in
\bbh\oplus\bbh^2\oplus\dots\oplus\bbh^2$ is any point such that
$q_j\neq 0$ for all $j$.  Then there is exactly one subspace from
the family for $G$ which contains  $q$. Such a point is clearly a
regular point of the cycle $G$.

We now consider a great circular curve $\mu\co  [0,\pi] \to
\bbp_{\bbr}^1$ which intersects the great circle defined by
$\ell(t)$ above, transversely (in one point). We then set
$$
\gamma(s) = (\mu(s), q_0,\ldots,q_0) \qquad\text{where}
\ q_0 = (1,0) \in \bbh^2.
$$
This closed curve intersects $G$ in exactly one point.  That point
is a regular point of $G$ and the intersection is transversal.
Thus by \cite[Lemma~9.13]{QQLLM} the cycle $G$ is not homologous
to zero and the proposition is proved.
\end{proof}
\medskip

To complete our analysis of the ring structure we need an explicit
representative of the generator $u$ of $\pi_2\cz_{\bbh}^{\ev}$.
Define the map
$$
\psi \co  S^2 = \bbp_{\bbc}^1(\bbh) \arr
\cz_{\bbh}^2(\bbp_{\bbc}(\bbh^2))
$$
by  $\psi = \widetilde{\psi} + \jj\widetilde{\psi}$ where
\begin{equation}
\label{eq:7.5}
 \widetilde{\psi}(\ell) = \ell_0 \oplus \ell - \ell_0 \oplus
\ell_0 \qquad\text{in}\ \bbh\oplus \bbh.
\end{equation}

To understand this map we  shall examine a basis for the homology
of $Q(\bbh^n)/\bbz_2$. Recall the following notation from Section
\ref{sec:6}:
$$
X^{2n-1} = \bbp_{\bbc}(\bbh^n)/\bbz_2  \aand  Y^{2n} =
Q(\bbh^n)/\bbz_2
$$
For $k\leq n$ there are   embeddings
$$
X^{2k-1}\subset Y^{2n}  \aand  Y^{2k} \subset  Y^{2n}
$$
where the second comes from the linear inclusion $\bbh^k\subset
\bbh^n$ and  $X^{2k-1}\subset Y^{2k}$ comes  from the 0--section of
${\mathcal O}(2)$. The analytic subsets  $Y^{2k}$ have oriented
regular sets  and define integral cycles which generate
$H_{4k}(Y^{2n};\, \bbz) = \bbz$. Each $X^{2k-1}$ is a smooth
nonorientable submanifold of the regular set of $Y^{2n}$. For
each $k$ consider the subspace $U^{2k+1}=\bbc\oplus \bbh^k \subset
\bbh^{n-k}\oplus\bbh^k$, let $Q(U^{2k+1})\subset Q(\bbh^n)$   be
the Thom space of  ${\mathcal O}_{\bbp(U^{2k+1})}(2)$, and set
$$
Z^{2k+1} = \pi\bigl(Q\bigl(U^{2k+1}\bigr)\bigr)
$$
where $\pi\co Q(\bbh^n)\to Y^{2n}$ is the projection.  Each
$Z^{2k+1}$ is an oriented analytic subvariety which defines an
integral cycle in $Y^{2n}$.

\begin{lemma}
\label{lem:7.5}
For each $k<n$ the class $[Z^{2k+1}]$  in $H_{4k+2}(Y^{2n};\,
\bbz) = \bbz_2$ is nonzero.
\end{lemma}
\begin{proof}
Note that $X^{2(n-k)-1}$ intersects  $Z^{2k+1}$ transversely in
exactly one point (in its regular set). The Lemma now follows from
\cite[Lemma~9.13]{QQLLM}.
\end{proof}
\medskip

\begin{proposition}
\label{prop:7.6}
Fix $m\geq 0$ and set
$$
F = \psi\# f \# \dots \# f \co  S^2\wedge  S^4\wedge \dots \wedge S^4
\arr \cz_{\bbh}^{2m+2}(\bbp_{\bbc}(\bbh^{2m+2})),
$$
where $f$ and $\psi$ are the maps from Proposition \ref{prop:7.1}
and \eqref{eq:7.5} respectively. Then $[F] = u\cdot x^m$ is the
generator of $\pi_{4m+2}\cz_{\bbh}^{\ev} = \bbz_2$.
\end{proposition}
\begin{proof}
Applying the homotopy equivalence $\csus=\csus_{{\mathcal O}(2)}$
of Proposition \ref{prop:6.1} gives a map
$$
\hF = \csus\circ F \co  S^{4m+2} \arr \cz^{2m+2}(Q(\bbh^{2m+2}))^\text{av}
$$
which splits as
$$
\hF = \tF +\widetilde{\jj} \tF
$$
where $\widetilde{\jj}$ is the real structure on $Q(\bbh^{2m+2})$.
 Proceeding in strict analogy with the proof of Proposition \ref{prop:7.4} we  are reduced
to showing that the cycle
$$\aligned
G &= \pi\biggl\{ \bigcup_{\ell, \la_1,\ldots,\la_m}
\csus\left(\ell_0\#\ell\#\la_1\#\dots\#\la_m  \right)\biggr\} \\
&= \pi\left\{ \csus\left( U^{2m+3} \right)\right\} = \pi\left(
Q( U^{2m+3})\right) = Z^{2m+3}
\endaligned
$$
is not 0 in $H_{4m+6}(Y^{2m+2};\,\bbz_2)$. This was proved in
Lemma \ref{lem:7.5}.
\end{proof}
\medskip

\noindent This completes the proof of Theorem \ref{thm:3.4}.


\section{Quaternionic projective varieties}
\label{sec:8}

The main theme of this paper is the study of spaces of
quaternionic cycles. Their structure turns out to be surprising
and rich.   However the geometry of quaternionic varieties
themselves is of  independent interest.  In this and subsequent
sections we will examine these varieties and show how our cycle
spaces provide invariants for their study.

\begin{definition}
\label{def:8.1}
 A \emph{quaternionic projective variety } is an
algebraic subvariety  $X \subset \bbp_{\bbc}(\bbh^n)$ which is
invariant under the quaternionic involution $\jj$.  A \emph{quaternionic
morphism} of quaternionic projective varieties is a
morphism which commutes with $\jj$.
\end{definition}

As we have seen, there are many quaternionic varieties.  The
abelian group they generate has the rich homotopy structure
determined above.  It is useful to look at some specific examples.

\begin{example}[Fermat varieties]
\label{ex:8.2}
Choose coordinates
$(q_1,\ldots,q_n)$ for $\bbh^n$ and write
\begin{equation}
\label{eq:8.1}
q_k = z_k + w_k\cdot j  \qquad\qquad\text{where} \ z_k, w_k
\in \bbc
\end{equation}
for all $k$.  Then the Fermat variety
$$
F(2m) = \biggl\{(z,w) \in \bbc^n \oplus\bbc^n\cdot j = \bbh^n
\,:\, \sum_k z_k^{2m} + w_k^{2m} = 0\biggr\}
$$
is a quaternionic variety for all $m\geq 1$. This includes the
K3--surface $F(4)\subset \bbp_{\bbc}(\bbh^2)$.
\end{example}

\begin{example}[Quaternionic divisors]
\label{ex:8.3}
Let $\Div_{2m} \cong
\bbp_{\bbc}^{\binom {2(n+m)}{2n} -1}$ be the space of divisors
of degree $2m$ on $\bbp_{\bbc}(\bbh^n)$.  Then $\jj$ induces a
linear antiholomorphic involution on $\Div_{2m}$ whose fixed-point
set is nonempty by Example \ref{ex:8.2}. Thus the subset
$\Div^{\bbh}_{2m} \subset\Div_{2m}$  of quaternionic divisors is a
real form
$$
\Div^{\bbh}_{2m} \cong \bbp_{\bbr}^{\binom {2(n+m)}{2n} -1}
\subset
  \bbp_{\bbc}^{\binom {2(n+m)}{2n} -1}.
$$
In general $\jj$ induces a real structure on $\Div_{\text{even}}$
and a quaternionic structure on $\Div_{\text{odd}}$.
\end{example}

The evenness of degree here is required by Theorem
\ref{thm:2.3}(ii). Here is an elementary proof.

\begin{proposition}
\label{prop:8.4}
Let $X \subset \bbp_{\bbc}(\bbh^n)$ be a quaternionic projective
variety of odd (complex) codimension. Then the degree of $X$ is
even.
\end{proposition}
\begin{proof}
 Let $\codim(X) = 2q-1$ and choose a quaternionic linear
subspace $V \subset \bbh^n$ with $\dim_{\bbh} = q$ such that
$\bbp(V)$ meets $X$ transversely at  regular points.  The
existence of such a  $V$ follows from the transitivity of $Sp_n$
on $\bbp_{\bbc}(\bbh^n)$ and Sard's Theorem for families (cf
\cite[Appendix~A]{QHL}).  Now $\bbp(V)\cap X$ is $\jj$--invariant, and
so $\deg(X) = \#(\bbp(V)\cap X)$ is even.
\end{proof}

\begin{example}[Quaternionic rational normal curves]
\label{ex:8.5}
Notice that the
mapping
$$
\bbp_{\bbc}(\bbh)  \arr \bbp_{\bbc}(\bbh^n)
$$
given by
$$
(z,w) \mapsto (Z,W) = (z^{2n-1}, z^{2n-2}w, \dots,
z^{n}w^{n-1}\,;\, w^{2n-1}, w^{2n-2}z, \dots, w^{n}z^{n-1})
$$
is a quaternionic morphism.  Its image is a $\jj$--invariant
rational normal curve. The moduli space of such curves is a real
form for the space of all rational normal curves in
$\bbp_{\bbc}(\bbh^n)$.
\end{example}

\begin{example}[Quaternionic Veronese and Segr\'e embeddings]
\label{ex:8.6}
More generally
one can check that there is a nonempty  subspace of \qa Veronese
embeddings
$$
\bbp_{\bbc}(\bbh^n) \arr
\bbp_{\bbc}(\text{Sym}_{\bbc}^{2k+1}\bbh^n)
$$
in any odd degree $2k+1$. This is also true of the Segr\'e
embeddings
$$
\bbp_{\bbc}(\bbh^{n_1})\times \dots \times \bbp_{\bbc}(\bbh^{n_k}) \arr
\bbp_{\bbc}(\bbh^{n_1}\otimes_{\bbc}\dots\otimes_{\bbc}\bbh^{n_k})$$
for all $k$ odd.
\end{example}

\begin{example}[General quaternionic curves]
\label{ex:8.7}
As we saw in Section~2 quaternionic
curves can be thought of as ``nonorientable''
 algebraic curves over $\bbc$.
\end{example}

\begin{theorem}[Intrinsic characterization of \qa projective manifolds]
\label{thm:8.8}
Let $X$ be a compact K\"ahler manifold with an antiholomorphic
involution $j\co X\to X$. Suppose there exists a positive holomorphic
line bundle $\pi\co L \to X$ which admits an anti-linear bundle map
$\widetilde j\co E\to E$ such that
$$
\bfig
\square<700,400>[L`L`X`X;\widetilde j`\pi`\pi`j]
\efig$$
commutes, and ${\widetilde j}^2 = -\Id$. Then there exists a $j$--equivariant
holomorphic embedding $\Phi\co X\hookrightarrow
\bbp_{\bbc}(\bbh^N)$ for some $N$.
\end{theorem}

\begin{proof}
Let $W_k$ denote the space of all holomorphic cross-sections of
the bundle $L^{\otimes k}$. By the fundamental theorem of Kodaira
for all $k$ sufficiently large, the mapping
$$
\Phi\co X \arr \bbp(W_k)^*\qquad \text{given by}\qquad \Phi(x) =
\ker\{\sigma \mapsto \sigma(x)\}
$$
is a well-defined projective embedding. For $k$ odd the bundle
$L^{\otimes k}$ is \qa and there is a \qa structure $\jj$ on $
W_k$ defined by setting $\jj(\sigma)\equiv {\widetilde
j}^{-1}\circ \sigma\circ j$. Note that $\Phi(jx) = \ker\{\sigma
\mapsto \sigma(jx)\} = \ker\{\sigma \mapsto {\widetilde
j}^{-1}\circ\sigma\circ j(x)\} = \{\jj(\sigma) : \sigma(x) =0\} =
\jj(\Phi(x))$, which proves the $j$--equivariance.
\end{proof}

This theorem motivates the definition of a \qa topological space
given below.


\section{Quaternionic algebraic cocycles and morphic
cohomology}
\label{sec:9}

We now want to consider families of quaternionic varieties $\pi\co \cf\to
X$ over a parameter space $X$.  Such families  generalize the concept of a
quaternionic vector bundle. We will  begin in the algebraic category and
adopt the viewpoint of algebraic cocyles developed in \cite{QFL}.  Much
of that theory carries over to the quaternionic case.

Recall the Chow monoid
$$\cc^q(\bbp_{\bbc}(\bbh^n)) = \coprod_{d\geq 0}
  \cc^q_d(\bbp_{\bbc}(\bbh^n))$$
where $\cc^q_d(\bbp_{\bbc}(\bbh^n))$ is the Chow variety of
effective algebraic cycles of degree $d$ and codimension $q$ in
$\bbp_{\bbc}(\bbh^n)$. The map $\jj$ induces an involution, also
denoted $\jj$, on each of these varieties. Let
$\cc^q_{\bbh}(\bbp_{\bbc}(\bbh^n))\subset
\cc^q(\bbp_{\bbc}(\bbh^n))$ denote the submonoid of $\jj$--fixed
cycles.

\begin{definition}
\label{def:9.1}
Let $X$ be a quaternionic variety (or more generally any real
variety with involution given by the action of $\Gal(\bbc/\bbr)$).
By a \emph{quaternionic  algebraic cocycle} on $X$ we mean a
$\jj$--equivariant morphism
$$ \varphi \co  X  \arr \cc^q(\bbp_{\bbc}(\bbh^n))$$
for some $n$. The set of all  quaternionic algebraic cocycles
 forms an abelian monoid
$$ \Mor_{\bbh}(X;\, \cc^q(\bbp_{\bbc}(\bbh^n)) $$
whose  group completion will be denoted by $ \Mor_{\bbh}(X;\,
\cz^q(\bbp_{\bbc}(\bbh^n)). $ Note that each cocycle $\varphi\in
\Mor_{\bbh}(X;\, \cc^q(\bbp_{\bbc}(\bbh^n))$
  gives rise to a mapping
\begin{equation}
\label{eq:9.1}
\widetilde{\varphi} \co  X/\bbz_2 \arr
\cz_{\bbh}^q(\bbp_{\bbc}(\bbh^n)) \qquad\qquad\text{where}
\ \widetilde{\varphi}([x]) = \varphi(x) + \jj\varphi(x).
\end{equation}
\end{definition}

\begin{example}[The fundamental class]
\label{ex:9.2}
Let $\varphi\co X\subset
\bbp_{\bbc}(\bbh^n)$ be a quaternionic variety.  This inclusion is
a \qa cocycle whose associated map
$$\widetilde{\varphi}\co X/\bbz_2\to
  \cz_{\bbh}^{2n-1}(\bbp_{\bbc}(\bbh^n))$$
is given by
$\widetilde{\varphi}([x]) = x + \jj x$. With respect to the
canonical splitting in Theorem \ref{thm:2.3}(ii) this fundamental
map can be viewed as
$$ \widetilde{\varphi}\co X/\bbz_2 \arr  \prod_{k=0}^{q} K(\bbz,
4k) \times \prod_{k=0}^{n} K(\bbz_2, 4k+1).
$$
Define total classes
$$
\iota = 1 + \iota_4 + \iota_8 + \iota_{12} +\cdots  \aand
\widetilde{\iota} = \widetilde{\iota}_2 + \widetilde{\iota}_6
+\widetilde{\iota}_{10} +\widetilde{\iota}_{14}+\cdots
$$
where $\iota_{4k}\in H^{4k}(K(\bbz, 4k);\, \bbz) = \bbz $ and $
\widetilde{\iota}_{4k+1} \in H^{4k+1}(K(\bbz_2, 4k+1);\, \bbz_2) =
\bbz_2 $ denote the fundamental classes.  Then associated to the
embedding $\varphi\co X\subset  \bbp_{\bbc}(\bbh^n)$ we have the
classes
\begin{equation}
\label{eq:9.2}
\widetilde{\varphi}^*  \iota \in H^{4*}(X/\bbz_2;\,\bbz) \aand
\widetilde{\varphi}^*\widetilde{\iota}\in
H^{4*+1}(X/\bbz_2;\,\bbz_2).
\end{equation}
When $X=\bbp_{\bbc}(\bbh^n)$, these classes are nonzero in every
dimension.

For a general \qa variety $X\subset \bbp_{\bbc}(\bbh^n) $ of
dimension $2q-1$ we can find quaternionic projection
$\bbp_{\bbc}(\bbh^n)-\bbp_{\bbc}(\bbh^{n-q-1}) \arr
\bbp_{\bbc}(\bbh^q)$ which restricts to give a \qa morphism $X
\arr \bbp_{\bbc}(\bbh^q)$.  The classes \ref{eq:9.2} for $X$ are
the pull-backs of those for $\bbp_{\bbc}(\bbh^q)$.
\end{example}

\begin{example}[The \qa Gauss map]
\label{ex:9.3}
Let $X\subset \bbp_{\bbc}(\bbh^n)$ be a
smooth quaternionic variety of codimension--$q$, and consider the
\qa morphism
$$
\g \co  X  \arr G^q_{\bbc}(\bbp_{\bbc}(\bbh^n))
\qquad\qquad\text{where} \ \g(x) = [T_xX].
$$
The associated characteristic map $ \widetilde{\g} \co  X/\bbz_2
\arr \cz^q_{\bbh}(\bbp_{\bbc}(\bbh^n)) $ can be rewritten in
terms of the canonical splitting in Theorem \ref{thm:2.3} as a
mapping
$$
\widetilde{\g} \co  X/\bbz_2  \arr \begin{cases} \prod_{j=0}^{r}
K(\bbz, 4j) \times
\prod_{j=0}^{r-1} K(\bbz_2, 4j+2)  &\text{if $q=2r$},\\
\prod_{j=0}^{r} K(\bbz, 4j) \times \prod_{j=0}^{r} K(\bbz_2, 4j+1)
&\text{if $q=2r+1$}.
\end{cases}
$$
Let $\iota$ and $\widetilde{\iota}$ be defined as in Example
\ref{ex:9.2} when $q=2r$ and let them be the obvious analogues
when $q=2r+1$. Then we can define the \emph{normal quaternionic
characteristic classes of $X$}:
\begin{equation}
\label{eq:9.3}
\widetilde{\g}^* (\iota) \in H^{4*}(X/\bbz_2;\, \bbz)
\qua\text{and}\qua
\widetilde{\g}^* (\widetilde{\iota}) \in \begin{cases}
H^{4*+2}(X/\bbz_2;\, \bbz_2)&\text{if $q=2r$}\\
H^{4*+1}(X/\bbz_2;\, \bbz_2)&\text{if $q=2r+1$}
\end{cases}
\end{equation}
As an example consider the Fermat variety $F(2) \subset
\bbp_{\bbc}(\bbh^n)$. Its  Gauss map
$$
\g \co  F(2)  \arr \bbp_{\bbc}^*(\bbh^n) =
G^1_{\bbc}(\bbp_{\bbc}(\bbh^n))
$$
is essentially the identity ($F(2)$ is self-dual). The associated
map
$$
\widetilde{\g} \co  F(2)/\bbz_2  \arr
\cz_{\bbh}^1(\bbp_{\bbc}(\bbh^n)) \cong K(\bbz_2,1)
$$
is easily seen to be nontrivial on $\pi_1$ and so
$\widetilde{\g}^* (\widetilde{\iota})\neq 0$.
\end{example}

\begin{example}
\label{ex:9.4}
Let $X\subset \bbp_{\bbc}(\bbh^n)$ be an irreducible hypersurface
of degree $d>1$.  The we can define
$$
\varphi\co  X \arr \cc_{d^2}^2(\bbp_{\bbc}(\bbh^n))
$$
by $\varphi(x) = X\bullet T_xX$ where  ``$\bullet$'' is the
intersection product (\cite{QFu}).
\end{example}

\begin{example}
\label{ex:9.5}
Consider the product \qa variety
$$
X = x_0\times \bbp_{\bbc}(\bbh) \amalg (\jj x_0)\times
\bbp_{\bbc}(\bbh) \cong \bbz_2\times S^2
$$
where $x_0$ is a point,  and the map $f\co X \to
G^2_{\bbc}(\bbp_{\bbc}(\bbh^2))$ given by
$$
f(x_0\times\ell) = \ell_0\oplus \ell \aand
 f(\jj x_0\times\ell) = \jj \ell_0\oplus \ell
\qquad\text{in}\ \bbh\oplus\bbh
$$
where $\ell_0 \subset \bbh\oplus \{0\} \subset \bbh^2$ is a fixed
complex line. Then the map
$$
\widetilde f \co  X/\bbz_2 = S^2  \arr \cz^2_{\bbh}
$$
represents the generator of $\pi_2\cz^2_{\bbh} = \bbz_2$ as we saw
in Proposition \ref{prop:7.6}.
\end{example}

\begin{example}
\label{ex:9.6}
Suppose $X\subset \bbp_{\bbc}(\bbh^2)$ is a \qa algebraic surface
of degree $2k$ which contains no \qa lines. Then there is a
well-defined continuous map
$$
\psi_X \co  S^4  \arr \cz^3_{\bbh}
$$
given by $\psi_X (p) = \pi^{-1}(p)\bullet X$ where
$\pi\co \bbp_{\bbc}(\bbh^2) \arr S^4$ is the Hopf fibration (See
Proposition \ref{prop:7.1}) and ``$\bullet$'' denotes intersection
product (cf \cite{QFu}). Now recall the isomorphism  $\tau\co
\cz^3_{\bbh}(\bbp_{\bbc}(\bbh^2)) \arr
\cz_{0}(\bbp_{\bbc}(\bbh^2)/\bbz_2) $ and set $\widetilde{\psi}_X
= \tau\circ\psi_X$.  Then if $\widetilde{\pi}\co
\cz_0(\bbp_{\bbc}(\bbh^2)/\bbz_2) \arr  \cz_0( S^4)$ is the
extension of the map $\bbp_{\bbc}(\bbh^2)/\bbz_2 \arr S^4$ then
$\widetilde{\pi}\circ\widetilde{\psi}_X = k\cdot\Id$.  Therefore
$$
[\psi] = k \in \pi_4\cz^3_{\bbh}.
$$
\end{example}


\section{Linear cocycles}
\label{sec:10}

The cocycles introduced  in Definition \ref{def:9.1} are
particularly interesting when $\varphi(x) $ is a linear subspace
for all $x$.

\begin{definition}
\label{def:10.1}
Let $X$ be as in Definition \ref{def:9.1}.  By an {\sl effective
quaternionic bundle of dimension $q$ } on $X$ we mean a
$\jj$--equivariant morphism $f\co X \to
G^q_{\bbc}(\bbp_{\bbc}(\bbh^N))$ for some $N$.
\end{definition}

Such a morphism corresponds to an algebraic vector bundle $E\to X$
which is generated by its global sections and which is equipped
with an anti-linear bundle map $\widetilde j\co E\to E$ which covers
$j\co X\to X$ and satisfies $\widetilde{j}^2=-\Id$.  These linear
cocycles form a submonoid under the algebraic join operation
($\cong$ direct sum in this case). The homotopy groups of its
group completion are interesting invariants of the variety. In
fact this group completion can be expanded to a generalized
equivariant cohomology theory attached to $X$. (See \cite{dSLii}.)

Note that to any \qa bundle $f\co X \to
G^q_{\bbc}(\bbp_{\bbc}(\bbh^N))$ there is  an associated  mapping
$$
\widetilde{f} \co  X/\bbz_2 \arr \cz_{\bbh}^q(\bbp_{\bbc}(\bbh^n))
$$
and we get classes $\widetilde{f}^*(\iota)$ and
$\widetilde{f}^*(\widetilde{\iota})$ as in \eqref{eq:9.3} above.
For the full theory of characteristic classes in this setting one
must consider the full equivariant theory. This is done in detail
in \cite{dSLi}


\section{\Qa spaces, \qa bundles and  $KH$--theory}
\label{sec:11}

The notions of  \qa vector bundles varieties can be generalized to
the topological category. Recall that a space with a \emph{real
structure} is a pair $(X,j)$ where $X$ is a topological space and
$j\co X\to X$ a continuous involution. The following notion was
introduced by  Johann Dupont \cite{Dui}.

\begin{definition}
\label{def:11.1}
A \emph{\qa vector bundle} over a {real space} $(X,j)$ is a
complex vector bundle $E\arr X$ together with an
$\bbc$--anti-linear bundle map $\widetilde{j}\co E\to E$ covering $j$
with $\widetilde{j}^2 = -1$.
\end{definition}

\noindent Such a pair $(E,\widetilde{j})$ with $\widetilde{j}^2=1$
is called a \emph{real bundle} (cf \cite{QA}).  Real bundles are
classified by $\bbz_2$--equivariant maps into the stabilized
Grassmannian with its standard real structure \cite{QQLLM}. The
corresponding statement holds for \qa bundles.

\begin{theorem}
\label{thm:11.2}
Let $X$ be a compact Hausdorff space with involution $j\co X\to X$.
Then the isomorphism classes of \qa vector bundles of complex rank
$q$ on   $X$ are in one-to-one correspondence with
$\bbz_2$--homotopy classes of $\bbz_2$--maps $X\to
G^{q}_{\bbc}(\bbp_{\bbc}(\bbh^{\infty}))$.
\end{theorem}
\begin{proof}
This is a direct adaptation of the standard arguments (cf
\cite{QM}).
\end{proof}
\begin{corollary}
\label{cor:11.3}
Let $E\to X$ be a \qa vector bundle classified by a
$j$--equivariant map $f\co X\to
G^{q}_{\bbc}(\bbp_{\bbc}(\bbh^{\infty}))$. Then the classes
$\widetilde{f}^*(\iota)$ and $\widetilde{f}^*(\widetilde{\iota})$
defined as in \eqref{eq:9.3} above, depend only on the isomorphism
class of $E$.
\end{corollary}

The Grothendieck group $KR(X)$ of real bundles on $X$ are the
basis of Atiyah's real $K$--theory \cite{QA}. The Grothendieck group
$KH(X)$ of \qa bundles on $X$ form an analogous \emph{\qa $K$--theory}
\cite{Dui}. However, this theory is not multiplicative. The tensor
product of two \qa bundles is  not \qa; it is real. However, the
combined theory $KR(X)\oplus  KH(X)$ has a product structure, and
interestingly there is an isomorphism
\begin{equation}
\label{eq:11.1}
KR(X\times \bbp_{\bbc}(\bbh)) \cong KR(X)\oplus KH(X)
\end{equation}
observed by Dupont \cite{Dui}. In a subsequent paper \cite{Duii}
Dupont asks whether there is an appropriate theory of
characteristic classes for \qa bundles and KH--theory. An answer to
the analogous question for real bundles and KR--theory was given by
dos Santos \cite{dSii}. The answer in the quaternionic case has
been recently given by dos Santos and Lima-Filho \cite{dSLi} who
continued this study of the space of \qa algebraic cycles.

\begin{definition}
\label{def:11.4}
A \emph{\qa space} is a triple $(X,j, \call)$ where $X$ is a
topological space, $j\co X\to X$ is an involution, and $\call \to X$
is a complex line bundle with \qa structure, that is, with a lifting
of $j$  to an anti-linear bundle map $\widetilde{j}\co \call \arr
\call$ such that $\widetilde{j}^2 = -1$. Note that $j$ must be a
free action.
\end{definition}

\begin{note}
 On a   \qa projective variety $X$ we take $\call = {\mathcal O}(1)$.
\end{note}

\begin{example}[Quaternionifications of a space]
\label{ex:11.5}
Any
real space $(X,j)$, with possibly trivial involution, gives rise
to  a \qa space $X_{\bbh}$ as follows.  Set $X_{\bbh} = \bbz_2
\times X$. Define $\jj\co X_{\bbh}\to X_{\bbh}$ by $\jj(0,x) =
(1,jx)$ and $\jj(1,x) = (0,jx)$. Set $\call = X_{\bbh}\times \bbc$
and define $\widetilde{\jj}\co \call\to\call $ by
$$
\widetilde{\jj}(0,x,z) = (1,jx,\bar z) \aand
\widetilde{\jj}(1,x,z) = (0,jx,-\bar z).
$$
This is the  \emph{trivial quaternionification} of $X$. There is a
natural bijection between complex bundles on $X$ and \qa bundles
on $X_{\bbh}$, and also between complex bundles on $X$ and real
bundles on $X_{\bbh}$. In particular we have $K(X) \cong
KR(X_{\bbh})\cong KH(X_{\bbh})$. In light of \eqref{eq:11.1} above
a more interesting quaternionification of $X$ is given by $X\times
\bbp_{\bbc}(\bbh)$ with $\call =\text{pr}_2^*{\mathcal O}(1)$.
\end{example}


\section{The equivariant homotopy type of
$\cz^q(\bbp_{\bbc}({\bbh^{\infty}}))$}
\label{sec:12}

There are two distinct real structures on projective space   (this
reflects the fact that the Brauer group of $\bbr$ is $\bbz_2$),
and they in turn induce real structures on the groups of algebraic
cycles. The first real structure, given by complex conjugation of
homogeneous coordinates, was studied in part one of this work
\cite{QQLLM}. The second, given by the quaternion involution on
$\bbp_{\bbc}(\bbh^n)$ (called the \emph{Brauer--Severi variety}), is
studied here. It is natural to ask: what is the full equivariant
homotopy type of the groups of algebraic cycles under the induced
involutions?

Recall from \cite{QL} that the nonequivariant homotopy type of
the group of cycles of codimension $q$  on $\bbp_{\bbc}^n$ is a
product of Eilenberg--MacLane spaces $K(\bbz,0)\times
K(\bbz,2)\times \dots \times K(\bbz,2q)$. In his thesis
\cite{dSii} Pedro dos Santos proved a beautiful, analogous result
for cycles under the first involution. He showed that there is a
$\bbz_2$--equivariant homotopy equivalence
$$
\cz^q(\bbp_{\bbc}(\bbc^n))\cong \prod_{k=0}^q K(\underline
\bbz, \bbr^{k,k})
$$
for any $n> q$, where $K(\underline \bbz,\bbr^{k,k}) $ denotes the
Eilenberg--MacLane space classifying $\bbz_2$--equivariant
cohomology indexed at the representation $\bbr^{k,k}$ ($=\bbc$ with
complex conjugation) with coefficients in the constant Mackey
functor $\underline \bbz$.

Very recently dos Santos and Lima-Filho \cite{dSLi} established
the corresponding result in the Brauer--Severi case. They prove
that there are $\bbz_2$--equivariant homotopy equivalences
$$
\cz^{2q-1}(\bbp_{\bbc}(\bbh^n))\cong \prod_{k=0}^q \text{Map}(
\bbp_{\bbc}(\bbh)_+,\,K(\underline \bbz, \bbr^{2k-1,\, 2k-1}))
$$
$$
\cz^{2q}(\bbp_{\bbc}(\bbh^n))\cong \prod_{k=0}^q \text{Map}(
\bbp_{\bbc}(\bbh)_+,\,K(\underline \bbz, \bbr^{2k,\, 2k})).
\leqno{\hbox{and}}$$
These spaces classify the (k,k)-equivariant cohomology of
$X\times\bbp_{\bbc}(\bbh)$.

%
%
%


\begin{thebibliography}

\bibitem{QA}  \textbf{M\,F Atiyah},
\emph{$K$--theory and reality},
Quart. J. Math. Oxford 17 (1966) 367--386
  \MR{0206940}

\bibitem{QBLLMM}
\textbf{C\,P Boyer}, \textbf{H\,B Lawson Jr}, \textbf{P Lima-Filho},
\textbf{B Mann}, \textbf{M-L Michelson},
\emph{Algebraic cycles and infinite loop spaces},
Invent. Math. 113 (1993) 373--388
  \MR{1228130}

\bibitem{QDT}
\textbf{A Dold}, \textbf{R Thom}, \emph{Quasifaserungen und unendliche
symmetrische produkte}, Ann. of Math. 67 (1956) 230--281
  \MR{0097062}

\bibitem{dSii}  \textbf{P\,F dos Santos},
\emph{Algebraic cycles on real varieties and $\bbz/2$--equivariant
homotopy theory}, Proc. London Math. Soc.  {86} (2003) 513--544
  \MR{1971161}

\bibitem{dSLi}
\textbf{P\,F dos Santos}, \textbf{P Lima-Filho},
\emph{Quaternionic algebraic cycles and reality},
Trans. Amer. Math. Soc. {356} (2004) 4701--4736
  \MR{2084395}

\bibitem{dSLii}
\textbf{P\,F dos Santos}, \textbf{P Lima-Filho},
\emph{Quaternionic $K$--theory for real varieties} (in preparation)

\bibitem{Dui}
\textbf{J Dupont},
\emph{Symplectic bundles and KR--theory},
Math.  Scand. {24} (1969) 27--30
  \MR{0254839}

\bibitem{Duii}
\textbf{J Dupont}, \emph{A note on characteristic classes for real
vector bundles}, preprint (1999)

\bibitem{QFL}
\textbf{E Friedlander}, \textbf{H\,B Lawson Jr},
\emph{A theory of algebraic cocycles},
Ann. of Math. {136} (1992) 361--428
  \MR{1185123}

\bibitem{QFu}
\textbf{W Fulton}, \emph{Intersection theory},  Ergebnisse series (3) 2,
Springer--Verlag,
New York (1984)
  \MR{0732620}

\bibitem{QHL}
\textbf{R Harvey}, \textbf{H\,B Lawson Jr},
\emph{On boundaries of complex analytic varieties I},
Ann. of Math. {102} (1975)  223--290
  \MR{0425173}

\bibitem{QLam}
\textbf{T-K Lam}, \emph{Spaces of real algebraic cycles and
homotopy theory}, PhD thesis, SUNY, Stony Brook (1990)

\bibitem{QL}
\textbf{H\,B Lawson Jr}, \emph{Algebraic cycles and homotopy theory},
Ann. of Math. {129} (1989) 253--291
  \MR{0986794}

\bibitem{QLLM}
\textbf{H\,B Lawson Jr}, \textbf{P Lima-Filho}, \textbf{M-L Michelsohn},
\emph{On equivariant algebraic suspension},
J. Algebraic Geom. {7} (1998) 627--650
  \MR{1642736}

\bibitem{QQLLM}
\textbf{H\,B Lawson Jr}, \textbf{P Lima-Filho}, \textbf{M-L Michelsohn},
\emph{Algebraic cycles and the classical groups I: Real cycles}, Topology {42}
(2003) 467--506
  \MR{1941445}

\bibitem{QLM}
\textbf{H\,B Lawson Jr}, \textbf{M-L Michelsohn},
\emph{Algebraic cycles, Bott periodicity, and the Chern characteristic
map}, from: ``The Mathematical Heritage of Hermann Weyl (Durham, NC,
1987)'', Proc. Sympos. Pure Math. 48, Amer. Math.
Soc., Providence, RI, USA (1988) 241--264
  \MR{0974339}

\bibitem{QLi}
\textbf{P Lima-Filho},
\emph{Lawson homology for quasiprojective varieties},
Compositio Math. {84} (1992) 1--23
  \MR{1183559}

\bibitem{QQLi}
\textbf{P Lima-Filho}, \emph{Completions and fibrations for
topological monoids}, Trans. Amer. Math. Soc.  {340} (1993)
127--146
  \MR{1134758}

\bibitem{QQQLi}
\textbf{P Lima-Filho}, \emph{On the equivariant homotopy of free
abelian groups on $G$--spaces and $G$--spectra}, Math. Z. {224}
(1997) 567--601
  \MR{1452050}

\bibitem{QMay}
\textbf{J\,P May}, \emph{${E}_\infty$ ring spaces and
${E}_\infty$ ring spectra}, Lecture Notes in Mathematics 577,
Springer--Verlag, New York, USA (1977)
\MR{0494077}


\bibitem{QM}
\textbf{J Milnor}, \textbf{J Stasheff}, \emph{Characteristic Classes},
Annals of Math. Studies 76, Princeton Univ. Press, Princeton, NJ, USA (1974)
  \MR{0440554}

\bibitem{QMo}
\textbf{J Mostovoy}, \emph{Quaternion flavored cycle spaces},
ICMS preprint, Edinburgh University (1975)

\end{thebibliography}
\end{document}